\documentclass[11pt]{article}
\usepackage[english]{babel}
\usepackage[utf8]{inputenc}
\usepackage{amsmath,amsfonts,amssymb,amsthm,mathrsfs,dsfont, bbm, mathtools}

\usepackage[mono=false]{libertine}
\useosf

\usepackage{wasysym}

\usepackage[dvipsnames]{xcolor}
\usepackage[a4paper,vmargin={3.5cm,3.5cm},hmargin={2.5cm,2.5cm}]{geometry}
\usepackage[font=sf, labelfont={sf,bf}, margin=1cm]{caption}

\usepackage{lscape}
\usepackage{pdflscape}
\usepackage{graphicx}

\usepackage[font=footnotesize,labelfont=bf]{caption}
\usepackage{subfig} 
\usepackage{epsfig}
\usepackage{latexsym}
\usepackage{stmaryrd}
\usepackage{ae,aecompl}

\usepackage[english]{babel}
 \usepackage[colorlinks=true]{hyperref}
\usepackage{pstricks}
\usepackage{enumerate}
\usepackage[normalem]{ulem}

\newtheorem{theorem}{Theorem}[]
\newtheorem{definition}{Definition}[]
\newtheorem{proposition}[]{Proposition}

\newtheorem{lemma}[]{Lemma}

\newtheorem{conjecture}[]{Conjecture}
\theoremstyle{definition}

\newtheorem*{remark}{Remark}

\renewcommand{\leq}{\leqslant}
\renewcommand{\geq}{\geqslant}

\newcommand{\Ind}[1]{\mathds{1}_{#1}}

\newcommand{\R}{\mathbb{R}}

\newcommand{\defeq}{:=} 

\makeatletter
\def\keywords{\xdef\@thefnmark{}\@footnotetext}
\makeatother

\newcommand{\K}{\mathtt K}

\begin{document}

\title{{\bf  A periodic Kingman model for the balance between mutation and selection.}}
\author{Camille Coron\thanks{Université Paris-Saclay, AgroParisTech, INRAE, UMR MIA Paris-Saclay, 91120, Palaiseau, France. Email: \texttt{camille.coron@universite-paris-saclay.fr}} \,  and \,
Olivier Hénard\thanks{Université Paris-Saclay, CNRS, Laboratoire de mathématiques d'Orsay, 91405, Orsay, France. Email: \texttt{olivier.henard@universite-paris-saclay.fr}}}

\maketitle

\begin{abstract}
We consider a periodic extension of the classical Kingman non-linear model \cite{Kingman_1978} for the balance between selection and mutation in a large population. In the original model, the fitness distribution of the population is modeled by a probability measure on the unit interval evolving through a simple dynamical system in discrete time: selection acts through size-biasing, and the mutation probability and distribution are kept fixed through time.
A natural extension of Kingman model is given by a periodic mutation environment; in this setting, we prove the convergence of the fitness distribution along subsequences and find an explicit criterion in terms of the Perron eigenvalue of an appropriately chosen matrix to decide whether an atom emerges at the largest fitness, a phenomenon usually called condensation.
\end{abstract}

\medskip
In population genetics, the equilibrium formed by the combination of selection and mutation is a major issue (\cite{Crow_Kimura_1971}, chapter 6.2). Indeed, while mutations have the potential to increase genetic diversity, selection concentrates the distribution of the genotypes around the fittest ones. Kingman house-of-cards model, introduced in \cite{Kingman_1978}, is precisely designed to study the balance between these two effects. In this model, individuals are characterized by a fitness $x\in[0,1]$ that simply determines their mean number of offspring. The main assumption of this model is that each new-born individual either inherits the fitness of its parent (chosen with probability proportional to its fitness), or mutates, in which case its fitness is drawn according to a fixed the mutation distribution, independently of its parent. This assumption is in sharp contrast with the small mutation assumption, yet it is sustained by some biological studies, see e.g. \cite{Hodgins2015} in which the authors show that Kingman house-of-cards model can explain the level of genetic diversity of several well-studied species. 

The two extreme cases are straightforward : Selection without mutation results in a Dirac mass at the largest fitness (all individuals share the same fitness at the limit), while mutation without selection results in a constant fitness distribution, the mutation distribution. In between, the effect of selection results in a limiting distribution skewed towards the largest fitness, plus, for a certain range a parameters, an atom at the largest fitness, an effect called condensation. Formally, Kingman house-of-cards model takes the form of a discrete dynamical system on the fitness distribution characterized by a mutation distribution and a mutation probability. The condensation phenomenon has also been studied in more general frameworks. In particular \cite{MaillerCourse} establishes conditions for condensation in probabilistic models of preferential attachment. In the deterministic side, in the case of a continuous-time parameter, the model, known as the replicator-mutator model, has been analysed in depth in \cite{cloez2024fast}, with also convergence rates; notice though the convergence is in Cesaro-mean only in the condensation case.
Slightly distinct models with a local mutation kernel (acting by convolution) that still display condensation have also been considered recently, see \cite{Giletal2017, Giletal2019}.

A generalization of Kingman house-of-cards model where the mutation probability is random is introduced in the series \cite{Yuan_2017, Yuan_2020, Yuan_2022}. There the author proves the convergence in distribution of the random fitness distribution, and derives criteria for condensation, both implicit and explicit (although arguably difficult to verify). More importantly for us, he also finds an important monotonicity property in Kingman's recursion that opens the road for a new proof of Kingman's theorem \cite{Henard_review} that is distinct from the ingenious but ad-hoc arguments of the original paper \cite{Kingman_1978}. 
 
In the present article we consider a different deterministic generalization of the model and concentrate on the case where the mutation distributions alternate in a periodic way, which leads again to a "solvable" model (in a sense, this is the main observation of this paper).
We prove the convergence of fitness distributions along sub-sequences and obtain an explicit criterion for condensation at the upper fitness, in terms of the Perron eigenvalue (aka spectral radius) of a matrix, as well as explicit formulae for the limiting fitness distributions. Relying on the properties of the matrix under consideration, we are able to phrase (perhaps suprisingly) our criterion in terms of the sign of a determinant, which is even simpler to check than the original criterion. A part is also devoted to the cumulative product of mean fitnesses which are particularly tractable quantities : we give closed expressions for their moment generating functions. Last, we point to further generalizations by phrasing a conjecture in the case of periodic mutation \emph{and} selection environments.

\section{Model and main result}

\subsection{Kingman original house-of-cards model}\label{sec-original-Kingman}

We first do a very quick review of Kingman house-of-cards model introduced in \cite{Kingman_1978} to motivate and give our results some context. We denote by $\mathcal P([0,1])$ the set of probability measures on $[0,1]$. Let $\beta$ be a real number in $(0,1)$, and let $p_0(dx), q(dx)\in\mathcal{P}([0,1])$ be respectively the initial fitness distribution, and the mutant fitness distribution. We set $\eta_q := \sup\{ x \in [0,1] : q([x, 1])>0 \}$ and $\eta_{0} := \sup\{ x \in [0,1] : p_0([x, 1])>0\}$ the respective maxima of the supports of $q(dx)$ and $p_0(dx)$, and assume that:
\begin{equation} 
\label{eq:ass-q} 
\eta_{0} \geq \eta_q >0.
\end{equation} 
We then define a sequence of probability measures on the unit interval $(p_n(dx))_{n\geq 0}$, where $p_n(dx)$ is the fitness distribution at time $n$,
 by Kingman's recursion:
 \begin{equation*}
  \left\{\begin{array}{ll}
 p_{n+1}(dx)  & = \beta q(dx) + (1- \beta) \frac{x p_{n}(dx)}{w_{n}}, \qquad \text{ for each integer }n \geq 0 \\
w_n   & = \int x p_n(dx) 
 \end{array}\right. \; 
\end{equation*}

Within assumption \eqref{eq:ass-q}, the part $\eta_q \leq \eta_0$ is without loss of generality since it holds after replacing $p_0(dx)$ by $p_1(dx)$, whereas the positivity assumption is to avoid uninteresting definition issues. The parameters in this model are the initial fitness distribution $p_0(dx)$, the mutation fitness distribution $q(dx)$ and the mutation probability $\beta\in (0,1)$, and one may gather the two parameters $\beta$ and $q(dx)$ in a sub-probability measure $\beta q(dx)$ that will be called the \emph{mutation environment} in the rest of the article.
We set, using the convention used all along the paper that all unspecified integrals will be on the unit interval,
$$\rho(z) = \sum_{n \geq 0} z^n  \int x^n \beta q(dx) = \int  \frac{\beta q(dx)}{1- zx}  \in (0,\infty], \text{ for each real } z \in [0,1/\eta_q],$$
and observe that the map $z \mapsto \rho(z)$ is continuous increasing on $[0,1/\eta_q]$.

\begin{theorem}[Kingman, 1987]
\label{main1} 
Either $\rho(1/\eta_0) \geq 1$, in which case the equation: 
$\rho(z) = 1$ with unknown $z \in [0, 1/\eta_0]$
has a unique solution denoted $z_c$; or $\rho(1/\eta_0) < 1$, in which case we set $z_c \defeq 1/\eta_0$.
Let also $\alpha= 1-\rho(z_c)$. We then have:
\begin{equation}
\label{conv-1d}
p_{n}(dx) \underset{n\to\infty}{\longrightarrow} \pi(dx)  \defeq \frac{\beta q(dx)}{1- z_c x}  + \alpha \delta_{\eta_0},
\end{equation}
for the weak convergence of probability measures on the interval $[0,\eta_0]$.
Furthermore, 
\begin{itemize}
\item[$(i)$] \emph{(non-condensation)} In case $\rho(1/\eta_0)\geq 1$, it holds $\alpha=0$ and the convergence also occurs in total variation.
\item[$(ii)$] \emph{(condensation)} In case $\rho(1/\eta_0)< 1$, it holds $\alpha>0$ and the convergence also occurs in total variation on any closed interval of $[0,\eta_0)$.
\end{itemize}
\end{theorem}

We stress that the limit $\pi(dx)$ of $(p_n(dx))_n$ is independent of $p_0(dx)$ in case $(i)$, and depends on $p_0(dx)$ only through the maximum $\eta_0$ of its support in case $(ii)$.
The most interesting case is $(ii)$ when an atom grows at the limit at $\eta_0$, meaning $\pi(\{\eta_0\}) \neq 0$, while neither $q(dx)$ nor $p_0(dx)$ had such an atom, $q(\{\eta_0\})=p_{0}(\{\eta_0\})=0$, an effect coined condensation; in this regard, we point out that the assumption $q(\{\eta_0\})\neq 0$ entails $\rho(1/\eta_0) = +\infty$, which implies that we lie in case $(i)$. For definiteness, we stress that in the rest of the article, the term condensation will refer to the fact we are in case $(ii)$. Let us point out that we can write in both cases 
\begin{equation}
\label{saturation-1d}
z_c  =  \sup\{ z \leq 1/\eta_0 : \rho(z) \leq 1 \}.
\end{equation}
The convergence in \eqref{conv-1d} entails the one of the mean fitnesses:
\begin{equation}
\label{saturation-1d-2}
w_n = \int x p_n(dx) \to \int x \pi(dx) = \frac{1-\beta}{z_c} \geq (1-\beta)\eta_0 \text{ as } n \to \infty
\end{equation}
with equality in the last inequality iff $z_c = 1/\eta_0$, that is $\rho(1/\eta_0) \leq 1$:  therefore, $(1-\beta)\eta_0$ may be viewed as a \emph{floor value} for the mean limiting fitness, with condensation occurring only if this value is reached. 
We end this review by noting that the model is rather generic: the choice of modeling the selection using size-bias with the identity function is indeed without loss of generality, since any other choice of function could be mapped to this one by an appropriate change of space (accordingly modifying the mutation environment).
The way we have phrased Theorem \ref{main1} is a bit distinct from the original statement in Kingman \cite{Kingman_1978} and closer to the statement of Yuan \cite{Yuan_2022}. Since Theorem \ref{main1} will appear as a special case of Theorem \ref{main2} that we prove in this article, we do not insist here on the precise connection.

\subsection{Kingman model with a periodic environment}
\label{periodic-sec}

Our aim in this article is to study the condensation phenomenon when the population experiences a periodic mutation environment. This means we consider $k$ mutation environments given by $k$ sub-probability measures on the unit interval $[0,1]$: $$(\beta_i q_i(dx))_{i \in \K}, \quad \text {  setting }  \quad \K = \{0,1, \ldots, k-1\},$$
where $\beta_i$ are real numbers in $(0,1)$ and $q_i(dx)$ are probability measures on the unit interval, for $i \in \K$.
For $\ell$ an integer, we set $[\ell]$ for $\ell \text{ mod } k$, the remainder of $\ell$ in the division by $k$, so that $[\ell]$ belongs to $\K $. The sequence $(p_n(dx))_n$ of fitness distributions is then defined inductively by:
\begin{equation}\label{eq:model} p_{n+1}(dx)  = \beta_{[n+1]} q_{[n+1]}(x) + (1-\beta_{[n+1]}) \frac{x p_{n}(dx)}{w_{n}}, \qquad \text{ for } n \geq 0.\end{equation}
We introduce as before the maximum of the support of the probability measures under consideration, namely we set $\eta_{0} \defeq \sup\{ x \in [0,1] : p_0([x, 1])>0\}$, $\eta_{q_i} \defeq \sup\{ x \in [0,1] : q_i([x, 1])>0 \}$ for $i\in\K$, and $\eta_{q}: =\max\{\eta_{q_i}, i \in \K\}$. As in the case with one environment $k=1$, we will work all along under the assumption that: 
\begin{equation}
\label{cond:eta}
\eta_0 \geq \eta_{q_i} >0 \; \text{, for each } \; i \in \K,
\end{equation}
in which the inequality between $\eta_0$ and $\eta_{q_i}$ is without any loss of generality (starting from $p_k(dx)$ instead of $p_0(dx)$ if needed).
The following notations are instrumental to state our main result.
Define 
\begin{equation*}
I_{\mu} = \bigcap_{j \in \K} \Big\{z \in \R_+:\sum_{n \geq 0} z^n \int x^n \beta_j q_j(dx)<\infty\Big\}.
\end{equation*}
Notice that $I_{\mu}$ is either $[0,1/\eta_q)$ or $[0,1/\eta_q]$.
Then we define for $i,j\in\K$ and $z\geq 0$ the following variant
of a moment generating function for the sub-probability measure $\beta_j q_j(dx)$: 
\begin{align}
\label{eq:mu-ij}
\mu_{j}^{i}(z) & :=  \sum_{n \geq 0 : [n]=i} z^n \int x^n \beta_j q_j(dx)  \in (0,\infty]  
\end{align}
and observe that $I_{\mu}$ also satisfies
$I_{\mu} =\{z \geq 0: \mu_j^i(z)<\infty, i,j \in \K\}$.
These quantities are in turn used to define a square matrix $A(z)$ with size $k$ parameterized by $z \in I_{\mu}$ as follows :
\begin{equation*}
A(z) = \bigg(\mu_{j}^{[i-j]}(z)\bigg)_{i,j \in \K}  \text{ where }
\mu^i_j(z) = \int \frac{(zx)^i}{1-(zx)^k} \beta_j q_j(dx),
\end{equation*}
where we stress that the integral representation on the right holds because $z \in I_{\mu}$. The matrix $A(z)$ is diagonal for $z=0$, with diagonal $(\beta_i)_{i \in \K}$, and has positive (finite) entries for $z\in I_{\mu}$. For $z \in I_{\mu}$, the Perron-Frobenius (or simply Perron) eigenvalue of $A(z)$ will be denoted by $\rho(A(z))$, see Appendix \ref{sec-appendix} for a statement of Perron theorem summing up the main properties of this ubiquitous object in our approach.
It may be the case that $1/\eta_0 \notin I_{\mu}$, but in this case we necessarily have $\lim_{z \to 1/\eta_0} \rho(A(z))= \infty$ (Lemma \ref{1eta0-def}) hence we shall set $\rho(A(1/\eta_0))=\rho(A(1/\eta_0)-)=\infty$.
In the next proposition  we formulate a dichotomy for the solution of an eigenvector equation that contains in germ the phase transition of Theorem \ref{main2}. We denote by $\mathbb R_+^\star=(0,\infty)$ the set of positive real numbers.

\begin{proposition}\label{prop:2cases}
 The matrix-valued function $A$ satisfies the alternative:
\begin{enumerate}
\item[(i)] Either $\rho(A(1/\eta_0)) \geq 1$, in which case the equation: 
\begin{equation}\label{eq:AX=X}A(z) X = X \end{equation}
with unknown $(X,z) \in (\R_+)^k \times I_{\mu}$ has a unique solution such that $X_0=1$; furthermore, this solution is such that $X \in \; (\mathbb R_+^\star)^k$.
\item[(ii)] Or $\rho(A(1/\eta_0)) < 1$, in which case $1/\eta_0 \in I_{\mu}$ and the equation :
$$A(1/\eta_0) X + \alpha \mathbf{1} = X $$
with unknown $(X,\alpha) \in \R^k \times \R$ has a unique solution such that $X_0=1$; furthermore, this solution is such that $(X,\alpha) \in \; (\mathbb R_+^\star)^k \times \; \mathbb R_+^\star$.
\end{enumerate}
\end{proposition}
\begin{remark}
We will show that in case $\rho(A(1/\eta_0))\geq 1$, Equation \eqref{eq:AX=X} still has a unique solution in the larger space $\R^k \times I_{\mu}$: this is the content of Proposition \ref{1-not-eigenvalue}, that we decided to single out in a separate statement because of its distinct proof that relies on the special form of the matrix $A$. 
\end{remark}

We can put the two items altogether as follows: defining the quantity $z_c$ by 
\begin{equation*}
z_c \defeq \sup\{ z \in [0,1/\eta_0] : \rho(A(z))\leq 1\},
\end{equation*}
(which is the periodic analog of \eqref{saturation-1d}), we shall observe that $z_{c}$ always belongs to $I_{\mu}$ (Lemma \ref{1eta0-def}) so that $A(z_c)$ is well defined, and Proposition \ref{prop:2cases} ensures
there is a unique vector $U$ with positive coordinates and first coordinate $U_0=1$ and a unique $\alpha  \geq 0$ satisfying :
\begin{equation}
\label{U-theta}
A(z_c)U+ \alpha \mathbf{1} = U.
\end{equation}
Notice that $\alpha= 0$ iff $\rho(A(1/\eta_0)) \geq 1$. With this piece of  notation at hand, we are able to state the following result, our main result in this article. Not only does it elucidate the convergence of the fitness distributions, but it also allows to decide for the apparition of an atom at the top of fitness distribution in terms of the just seen dichotomy for the matrix $A(1/\eta_0)$.

\begin{theorem}
\label{main2} 
We have, for $\alpha$ and $U$ defined in \eqref{U-theta}: 
$$
p_{kn+i}(dx) \underset{n\to\infty}{\longrightarrow} \pi_i(dx)  \defeq \sum_{j=0}^{k-1} \frac{U_{[i-j]}}{U_{i}} \frac{(z_c x)^j \beta_{[i-j]} q_{[i-j]}(dx)}{1- (z_c x)^k}  + \alpha\frac{U_{0}}{U_{i}} \delta_{\eta_0},   \qquad \text{for } i \in \K,
$$
for the weak convergence of probability measures on the interval $[0,\eta_0]$.
Furthermore, 
\begin{itemize}
\item[$(i)$] \emph{(non-condensation)} In case $\rho(A(1/\eta_0))\geq 1$, it holds $\alpha=0$ and the convergence also occurs in total variation on $[0,\eta_0]$.
\item[$(ii)$] \emph{(condensation)} In case $\rho(A(1/\eta_0))< 1$, it holds $\alpha>0$ and the convergence also occurs in total variation on any closed interval of $[0,\eta_0)$.
\end{itemize}
\end{theorem}

Before commenting on our main result, we point out that standard linear algebra considerations allow to be fairly explicit about $U$, namely we are able to give a closed formula in terms of the determinant of a matrix related to $A(z_c)$. 
If $\mathbf{I}_k$ stands for the identity matrix of dimension $k$, this matrix is defined for $z \in I_{\mu}$ by:
\begin{equation}
    \label{def:rhoj}
    B(z)  =A(z)-\mathbf{I}_k -\frac{1}{k}\Big(\rho_j(z)-1\Big)_{i,j\in \K} 
    \text{ where } \rho_j(z)= \int \frac{\beta_jq_j(dx)}{1-zx}. 
\end{equation}

\begin{proposition}
\label{closed-formula}
If $N_j$ denotes for every $j \in \K$ the $(j,j)$-minor of the matrix $B(z_c)$ (the determinant of the matrix with $j$-th lines and columns both removed), one has $N_0\neq0$ and, for each $j \in \K$,
$$U_j=\frac{N_j}{N_0}   \text{ and also } \alpha=\frac{1}{k}\sum_{j=0}^{k-1} \frac{N_j}{N_0}\left(1-\rho_j(z_c)\right).$$
\end{proposition}
Notice that, since only quotients of coordinates of $U$ appear in Theorem \ref{main2}, one can simply replace $U_i$ by $N_i$ in the expression in Theorem. Last but not least, relying on the properties of the matrix $A(z)$, we are able to establish  that there are no real eigenvalues $\geq 1$ apart from the Perron one, which results in an even simpler criterion for the condensation phase. 

\begin{proposition}
\label{equivalence2conditions}
Let $z \in I_{\mu}$. The two quantities  
$$ 1-\rho(A(z))  \quad \text{ and }\quad \det(\mathbf{I}_k-A(z))    $$
have the same sign, in the sense they are both negative, null or positive.
\end{proposition}

\begin{remark}
The interest of this proposition lies in its application to the condensation criterion; for this one needs to evaluate both sides at $z=1/\eta_0$. Provided $1/\eta_0 \in I_{\mu}$, condensation then occurs if only if and only if $\det(\mathbf{I}_k-A(1/\eta_0))>0$. In case $1/\eta_0 \notin I_{\mu}$, some coefficients of $A(1/\eta_0)$ are infinite and we will prove in Lemma \ref{1eta0-def} that there is no condensation.
\end{remark}

First, it is easy to check that specializing Theorem \ref{main2} to the $k=1$ case allows to recover Theorem \ref{main1}.
Theorem \ref{main2} shows that determining the emergence of condensation ($\alpha>0$) reduces to the computation of the Perron eigenvalue of the matrix $A(1/\eta_0)$ which in turn boils down to that of the determinant of the matrix $\mathbf{I}_k-A(1/\eta_0)$ according to Proposition \ref{equivalence2conditions}. Notice also that the quantity $\alpha$ keeps a clear interpretation as the mass of the atom at $\eta_0$ for the probability measure $\pi_0(dx)$. Theorem \ref{main2} (together with Proposition \ref{closed-formula}) also gives an explicit formula for the limiting fitness distribution along periodic sub-sequences. We point out that in the definition of $\pi_i(dx)$, the division by $U_{i}$ is legitimate since all coordinates of $U$ are positive.

In the most interesting case where $\eta_0$ is an atom of none of the $q_{i}(dx)$, $i \in \K$, then: either $\eta_0$ is an atom of none of the $\pi_{i}(dx)$ and each of the $\pi_{i}(dx)$ is absolutely continuous with respect to $\sum_{i} q_i$, case $(i)$,  
or $\eta_0$ is an atom of each of the $\pi_{i}(dx)$ and none of the $\pi_{i}(dx)$ is absolutely continuous with respect to $\sum_{i} q_i$, case $(ii)$. Last, if $\eta_0$ is an atom of $q_{i}(dx)$ for some $i$, we land necessarily in case $(i)$.


The limiting fitness distributions $\pi_i(dx)$ depend upon the initial distribution $p_0(dx)$ only via $\eta_0$ the maximum of its support, hence the limiting fitness distributions are invariant under cyclic permutations of the environments (after relabeling), but not in general by any permutation. In this same vein, we point to a direct proof (Lemma \ref{invariance-rotation}) that the quantity $\text{det}(\mathbf{I}_k-A(z))$ and $\rho(A(z))$ are invariant under cyclic permutations of the environments.

Concerning assumption \eqref{cond:eta}, we believe most of our results would still hold under the weaker condition ($\eta_0 >0$ and $\eta_0 \geq \eta_{q_i}$, for each $i \in \K)$, but this setting would require some additional work that is left to the interested reader.

Last, we point out that in the periodic setting, the choice of the identity function for the size-bias is no more without loss of generality (as was mentioned at the end of Section \ref{sec-original-Kingman} for the original $k=1$ model), and this prompts to also look for a model where the selection is itself periodic: we touch upon this subject through the formulation of Conjecture \ref{conjecture}.

\subsection{Comparison with the criterion for a single environment}

We conclude by commenting on the comparison between the criterion for periodic mutation environments and the original criterion by Kingman for a single environment. To that aim, we compare 
the Perron eigenvalue $\rho(A(z))$ of the matrix $A(z)$ to the values $\rho_j(z)$ defined in \eqref{def:rhoj}.

\begin{proposition}
\label{1-equivalent}
Let $z \in I_{\mu}$.
\begin{enumerate}
\item It holds:
\begin{equation}
\label{bounds-rho}
    \min\big\{ \rho_j(z) ,  j \in \K \big\} \leq  \rho(A(z)) \leq \max \big\{\rho_j(z) ,  j \in \K\big\}.
\end{equation}
\item  If the Perron eigenvector $R$ of the matrix $A(z)$ is normalized so that $\sum_{i \in K} R_i=1$, the two quantities:
\begin{equation*}
    \rho(A(z)) - 1 \quad \text{ and } \quad \sum_{i \in K} R_i \rho_i(z)- 1
\end{equation*}
have the same sign, meaning they are both negative, null or positive.
\end{enumerate}
\end{proposition}

\begin{remark}
Setting $z= 1/\eta_0$ leads to the following remarks.
Regarding the first point, if none of the environments leads to condensation in Kingman model with a single environment (i.e. $\rho_i(1/\eta_0) \geq 1$ for every $i \in \K$) then no condensation occurs in the model with periodic environment, and conversely if each environment leads to condensation in Kingman model with a single environment (i.e, $\rho_i(1/\eta_0) < 1$ for every $i \in \K$), then condensation will occur in the model with periodic environment. Also, taking $k$ equal environments $\beta q (dx)$, we have that the 3 terms in the inequality \eqref{bounds-rho} are equal, and recover Kingman criterion for a single environment.
Regarding the second point, taking $z = 1/\eta_0$, the periodic environment and the single environment given by the sub-probability measure $\sum_{i=0}^{k-1} R_i \beta_i q_i(dx)$ have the same criterion for condensation.
\end{remark}

\begin{proof}
For the first point: Using that $A(z)$ and its transpose have the same eigenvalues, hence the same Perron eigenvalue, and applying the bound given in Proposition \ref{prop-Perron} point (ii), one gets
$$\min_{j \in \K}\sum _{i} A_{ij}(z) \leq \rho(A(z)) \leq \max_{j \in \K} \sum _{i}A_{ij}(z),$$
but the sum of the coefficients of the $j$-th column of $A(z)$ is simply $\rho_j(z)$:
\begin{equation}
\label{sum:columns}
\sum_{i \in \K} A_{ij}(z) = \sum_{i \in \K} \mu_{j}^{[i-j]}(z) = \sum_{i \in \K} \mu_{j}^{i}(z)  = \sum_{i \in \K}  \int \frac{(zx)^i}{1-(zx)^k} \beta_j q_j(dx) = \rho_j(z).
\end{equation}
For the second point: if $\rho(A(z))> 1$, the Perron eigenvector $R$ of $A(z)$ satisfies: $A(z)R= \rho(A(z)) R > R$ coordinate-wise, 
hence, using \eqref{sum:columns} and the assumption that $R$ is normalized so that the sum of its coefficients is 1, we get $$\sum_{j} \rho_{j}(z) R_j= \sum_{i,j} A_{i,j}(z) R_j  = \rho(A(z)) \sum_{j}  R_j  >  1.$$
The cases where $\rho(A(z))= 1$ and $\rho(A(z))< 1$ are treated similarly. 
\end{proof}

This is the outline of the paper: Section \ref{sec:proofs} contains the proof of our main result, Theorem \ref{main2},  together with the results surrounding it. Following the line of reasoning developed in \cite{Yuan_2020}, we first prove (Proposition \ref{prop:convergence-fitness}) that for each $i\in \K$, the sequence of mean fitnesses $(\int x p_{kn+i}(dx))_{n\geq 0}$ converges. Next we prove (Proposition \ref{prop-cvgce-p}) that this convergence implies the convergence (in an appropriate sense) for each $i \in \K$, of the entire fitness distributions $(p_{kn+i}(dx))_{n\geq 0}$. We finally identify the limit as a solution of a fixed point equation (Proposition \ref{prop-equation}). The proof of Proposition \ref{prop:2cases} then gives the missing uniqueness statement to close the proof of Theorem \ref{main2}. Last, algebraic manipulations allow to give the closed formulae enounced in Proposition  
\ref{closed-formula} and conclude the section.
Section \ref{sec-additional} contains some additional elements: first, the alternative criterion in terms of the determinant of $\mathbf{I}_k-A(1/\eta_0)$, Proposition \ref{equivalence2conditions}, is proved.
Then some generating functions of cumulative products of mean fitnesses are computed, in the vein of the original approach by Kingman.  Last, a conjecture is phrased for an equivalent of Theorem \ref{main2} in the case where the mutation and selection environments are periodic.  
An appendix gathers a few results on matrices, and notably a key result, Lemma \ref{1-not-eigenvalue}, on the non-existence of eigenvectors with coordinates of different signs associated with real eigenvalues $\geq 1$ for a special class of matrices with positive coefficients that notably contains the matrices $A(z)$ for $z>0$.

\section{Proofs}\label{sec:proofs}

\subsection{Convergence of the sequence of fitness distributions along subsequences}

In this section, we prove that, for each $i \in \K$, the subsequence $(p_{kn+i}(dx))_{n}$ converges (in an appropriate sense) as $n \to \infty$ towards
$\pi_i(dx)$ with mean $\bar{w}_i = \int x \pi_{i}(dx)$. Our proof consists in first assessing the convergence of the mean fitness subsequences $(w_{kn+i}(dx))_{n}$ and then proving this convergence implies the convergence of the fitness distributions along subsequences.
To this aim we introduce after \cite{Yuan_2020} a partial order on the set of probability measures on the unit interval that has the particularity to be preserved through Kingman's recursion.
We parallel the arguments of \cite{Henard_review} that deals with the simpler case with one environment ($k=1$) only.

\begin{definition}
For $p, p' \in \mathcal{P}([0,1])$  and $\eta\in(0,1]$, we define a partial order between $p$ and $p'$ by: $$p\geq_{\eta^{-}}p'\quad\Leftrightarrow\quad p(A)\geq p'(A) \quad\text{ for every measurable set } A \text{ of } [0,\eta).$$
\end{definition}

Let us denote by $\Theta$ the Kingman's recursion operator with mutation probability $\beta$ and mutant fitness distribution $q\in\mathcal{P}([0,1])$: for any probability measure $p$ on $[0,1]$, 
$$\Theta(p)(dx)=\beta q(dx)+(1-\beta) \frac{x p(dx)}{\int x p(dx)}.$$
The following proposition collects basic results established in \cite{Yuan_2020} and gathered in this form in \cite{Henard_review}. 
We shall need the following notion of truncation of a measure: for $\varepsilon\in(0,\eta_0]$, and $\mu \in \mathcal P([0,\eta_0])$, we define $\mathcal{R}_{\eta_0-\varepsilon}$ the truncation of $\mu$ at $\eta_0-\varepsilon$ by: $\mathcal{R}_{\eta_0-\varepsilon}(\mu)(dx)=\mu|_{[0,\eta_0-\varepsilon)}(dx)+\mu([\eta_0-\varepsilon,\eta_0])\delta_{\eta_0-\varepsilon}.$

\begin{proposition}\label{lem:order}  Let $\beta\in(0,1)$, $q\in\mathcal{P}([0,1])$, and $p,p'\in\mathcal{P}([0,1])$ be such that their supports satisfy: $\eta:= \max(\emph{Supp}(p))\geq\max(\emph{Supp}(p'))$, then:
\begin{itemize}
\item[$(i)$] (point $(v)$ Lemma 5 of \cite{Henard_review}) If $p'\geq_{\eta^{-}}p$ then $\int xp'(dx)\leq \int xp(dx)$.
\item[$(ii)$] (Lemma 6 of \cite{Henard_review}) If $p'\geq_{\eta^{-}}p$ then $\Theta(p')\geq_{\eta^{-}} \Theta(p)$.
\item[$(iii)$] (Corollary 4 of \cite{Henard_review}) Let $\varepsilon\in(0,\eta_0]$. The sequence $(p'_{n}(dx))_{n\geq 0}$ defined by  
$p'_0=\mathcal{R}_{\eta_0-\varepsilon}(p_0)$  and
\begin{equation*} 
p'_{n+1}(dx)  = \mathcal{R}_{\eta_0-\varepsilon}(\beta_{[n+1]}q_{[n+1]} + (1-\beta_{[n+1]}) \frac{x p'_{n}(dx)}{\int x p'_n(dx)},  \qquad \text{ for } n \geq 0,
\end{equation*}
satisfies:
$$p'_n\geq_{(\eta_0-\varepsilon)^-}p_n \quad\text{ for all $n\geq 0$},$$
where $(p_n(dx))_{n\geq 0}$ satisfies the periodic Kingman model \eqref{eq:model} and $\eta_0= \max(\emph{Supp}(p_0))$.
\end{itemize}
\end{proposition}

Solving Kingman recursion as given in Equation \eqref{eq:model} gives rise to the following expansion.
\begin{lemma}\label{lem:decomposition}
For each integer $n\geq 1$, the fitness distribution $p_n(dx)$ can be decomposed in two parts, that will be denoted $a_n(dx)$ and $b_n(dx)$ in the following\footnote{
we use the convention that products over empty set are 1 (here for the term associated with $\ell=0$)}:
\begin{align}
\label{eq-a_nb_n} 
p_n(dx)
& = 
\underbrace{\sum_{l=0}^{n-1}\left(\prod_{ 0 \leq i \leq l-1}\frac{(1-\beta_{[n-i]})x}{w_{n-i-1}}\right)\beta_{[n-l]}q_{[n-l]}(dx)}_{\defeq a_n(dx)}+ \underbrace{\left(\prod_{i=0}^{n-1}\frac{(1-\beta_{[n-i]})x}{w_{n-i-1}}\right)p_0(dx)}_{\defeq b_n(dx)}.\end{align}
\end{lemma}

\begin{proof}
    This result is proved by induction on $n$. It holds for $n=1$, and if true for a given integer $n\geq 1$, then writing
    \begin{align*}p_{n+1}(dx)&=\beta_{[n+1]}q_{[n+1]}(dx)+(1-\beta_{[n+1]})\frac{xp_n(dx)}{w_n}\\&=\beta_{[n+1]}q_{[n+1]}(dx)+\sum_{l=0}^{n-1}\left(\frac{(1-\beta_{[n+1]})x}{w_n}\prod_{i=0}^{l-1}\frac{(1-\beta_{[n-i]})x}{w_{n-i-1}}\right)\beta_{[n-l]}q_{[n-l](dx)}\\&\phantom{\beta_{[n+1]}q_{[n+1]}(dx)}+\left(\prod_{i=0}^{n}\frac{(1-\beta_{[n+1-i]})x}{w_{n+1-i-1}}\right)p_0(dx)\\&=\sum_{l=0}^{n}\left(\prod_{i=0}^{l-1}\frac{(1-\beta_{[n+1-i]})x}{w_{n+1-i-1}}\right)\beta_{[n+1-l]}q_{[n+1-l]}(dx)+\left(\prod_{i=0}^{n}\frac{(1-\beta_{[n+1-i]})x}{w_{n+1-i-1}}\right)p_0(dx)\end{align*} gives the result at $ n+1$.
\end{proof}

The following proposition is the key to the proof of the convergence of the fitness distribution. 
\begin{proposition}\label{prop:convergence-fitness}
\begin{enumerate}
    \item For each $i\in \K$, the sequence of mean fitnesses $(w_{kn+i})_{n\geq0}$ converges towards a quantity denoted $\bar{w}_{i}$ depending on $p_0$ only through $\eta_0$ the maximum of its support. 
    \item Furthermore, we have the following bound:
\begin{equation}
\label{ineq:w}
\eta_0^k \prod_{i\in\K}\frac{1-\beta_i}{\bar{w}_i } \leq 1.
\end{equation}
\end{enumerate}
\end{proposition}

The proof of point 1 follows the strategy developed in \cite{Henard_review}, that itself takes its root in \cite{Yuan_2022}.
\begin{proof}

The proof decomposes in four steps, labelled from $(i)$ to $(iv)$ depending on the initial condition $p_0$.

$(i)$ Let $(p_n^{\bullet}(dx))_{n\geq 0}$ be the sequence of fitness distributions defined by Kingman recursion \eqref{eq:model} with initial condition $p_0^{\bullet}=\delta_{\eta_0}$. By definition, $p_0^{\bullet}\leq_{\eta_0^-}p_k^{\bullet}$.
Then, from point $(ii)$ of Proposition \ref{lem:order}, applying one step of Kingman recursion with the environment $\beta_1q_1(dx)$ gives that $p_1^{\bullet}\leq_{\eta_0^-}p_{k+1}^{\bullet}$. A finite induction based again on point $(ii)$ (with the appropriate mutation environment) 
then gives  
$p_n^{\bullet}\leq_{\eta_0^-}p_{k+n}^{\bullet}$, and in particular, $$p_{i}^{\bullet}\leq_{\eta_0^-}p_{k+i}^{\bullet}\leq_{\eta_0^-}p_{2k+i}^{\bullet}\leq_{\eta_0^-}...\leq_{\eta_0^-}p_{nk+i}^{\bullet}$$
holds for any $i\in \K$ and any $n\geq 0$. Therefore the sequence of distributions $(p_{kn+i}(dx))_{n\geq 0}$ is non-decreasing for the order $\leq_{\eta_0^-}$, and then from point (i) of Proposition \ref{lem:order}, the sequence of mean fitnesses $(w_{kn+i}^{\bullet})_{n\geq0}$ is non-increasing, hence admits a limit that we denote by $\bar{w}^{\bullet}_i$.

Before proceeding with the rest of the proof of point 1, let us observe that, from Lemma \ref{lem:decomposition},
\begin{align}\label{eq:p(eta_0)} 1\geq p_{kn}^{\bullet}(\{\eta_0\})&\geq  
\prod_{j=0}^{kn-1}\frac{(1-\beta_{[j+1]})\eta_0}{w_{j}^{\bullet}}
=\prod_{m=1}^n\frac{ \prod_{l\in\K}((1-\beta_l)\eta_0)}{w^{\bullet}_{k(m-1)}...w_{km-1}^{\bullet}}\end{align}
and we just proved the sequence $(w^{\bullet}_{k(m-1)}...w_{km-1}^{\bullet})_{m\geq1}$ converges to $\prod_{l\in \K}\bar{w}_l^{\bullet}$; the upper bound by 1 in \eqref{eq:p(eta_0)} then implies
$\prod_{l\in \K}\bar{w}^{\bullet}_l \geq \prod_{l\in\K}((1-\beta_l)\eta_0)$, and this is point 2 of the Proposition with $\bar{w}_i$ replaced by $\bar{w}_i^{\bullet}$.

Now, we have proved that the sequence $(w^{\bullet}_{k(m-1)}...w_{km-1}^{\bullet})_{m\geq1}$ is non-increasing, and, by point 2, that its limit is greater than or equal to $\prod_{l\in\K}((1-\beta_l)\eta_0)$. This gives $w^{\bullet}_{k(m-1)}...w^{\bullet}_{km-1}\geq \prod_{l\in\K}((1-\beta_l)\eta_0)$ for any $m\geq 0$ hence the sequence on the RHS of \eqref{eq:p(eta_0)} is non-increasing in $n$ and  converges to a constant 
$$C:= \lim_{n}  \prod_{j=0}^{kn-1}\frac{(1-\beta_{[j+1]})\eta_0}{w_{j}^{\bullet}}\geq 0$$

We now distinguish $3$ cases, according to whether $C=0$ or ($C>0$ and $p_0(\{\eta_0\})>0$), or ($C>0$ and $p_0(\{\eta_0\})=0$): in each case, we prove that the sequence $(w_{kn+i})_n$ converges towards $\bar{w}_i^{\bullet}$, which will conclude the proof of point 1 and 2 (with $\bar{w}_i=\bar{w}_i^{\bullet}$).

$(ii)$ Case $C=0$. Using the notation $a_n^{\bullet}(dx)$ and $b_n^{\bullet}(dx)$ for the measures $a_n(dx)$ and $b_n(dx)$ defined in Equation \eqref{eq-a_nb_n} associated with $p_0=\delta_{\eta_0}$, we get for each $i\in\K$, 
$$\int b_{kn+i}^{\bullet}(dx) = \prod_{j=0}^{kn-1}\frac{(1-\beta_{[j+1]})\eta_0}{w_{j}^{\bullet}}
\prod_{u=1}^{i}\frac{(1-\beta_{u})\eta_0}{w^{\bullet}_{kn-1+u}}\rightarrow_{n\rightarrow\infty} C \prod_{u=1}^{i}\frac{(1-\beta_{u})\eta_0}{\bar{w}^{\bullet}_{u-1}}=0,$$
hence $\int  x b_{n}^{\bullet}(dx) \leq \int b_{n}^{\bullet}(dx) \to 0$. Now, $p_0\geq_{\eta_0^-} \delta_{\eta_0}$ therefore from point $(ii)$ of Proposition  \ref{lem:order}, $w_{n}\leq w_n^{\bullet}$ for all $n$ and this gives the first (by definition) and third inequality in the following chain:
$$\int x p_n^{\bullet}(dx) \geq \int x p_n(dx) \geq \int x a_n(dx) \geq \int x a_n^{\bullet}(dx) = \int x p_n^{\bullet}(dx) - \int x b_n^{\bullet}(dx).$$
This altogether implies that for each $i \in \K$, the sequence of mean fitnesses 
$(w_{kn+i})_{n \geq 0} = (\int x p_{kn+i}(dx))_{n \geq 0}$  
has the same limit $\bar{w}^{\bullet}_i$ as $n \to \infty$ as the sequence $(w_{kn+i}^{\bullet})_{n \geq 0}$.

$(iii)$ Case $C>0$ and $p_0(\{\eta_0\})>0$. 
Let $\gamma\in(0,1)$ and set $N_n=\text{Card}(\{l\leq n-1 : w_l\leq (1-\gamma) w_l^{\bullet}\}$. We start with the following lower bound: 
\begin{align*}p_{kn} (\{\eta_0\})& \geq
\prod_{j=0}^{kn-1}\frac{(1-\beta_{[j+1]})\eta_0}{w_{j}^{\bullet}} \left(\frac{1}{1-\gamma}\right)^{N_{kn}} p_0(\{\eta_0\}).\end{align*}

Since $\lim\prod_{i=0}^{kn-1}\frac{(1-\beta_{[n-i]})\eta_0}{w_{n-i-1}^{\bullet}}=C>0$, $(N_{kn})_{n \geq 0}$ (hence also $(N_{n})_{n \geq 0}$) must be a bounded sequence, since otherwise one would have  $p_{kn}(\{\eta_0\}) \to \infty$. 
Since this holds for each $\gamma\in(0,1)$, we get  $\liminf_n w_{kn+i} \geq \liminf w_{kn+i}^{\bullet} =\bar{w}_{i}^{\bullet}$; on the other hand
$w_n \leq w_n^{\bullet}$ ensures that the complementary bound  $\limsup_n w_{kn+i} \leq \bar{w}_{i}^{\bullet}$ holds.

$(iv)$ If $C>0$ and $p_0(\{\eta_0\})=0$, then we consider the rabotted version of Kingman recursion, namely we set $(p^{(\varepsilon)}_{kn+i}(dx))_{n\geq0}$ the sequence of fitness distribution satisfying Kingman recursion starting from $\mathcal{R}_{\eta_0-\varepsilon}(p_0)$ with periodic mutation environments $(\mathcal{R}_{\eta_0-\varepsilon}(\beta_{[n]}q_{[n]}))_{n\geq1}$, and we observe it falls into one of the previous cases since $\mathcal{R}_{\eta_0-\varepsilon}(p_0)(\{\eta_0-\varepsilon\})>0$ by definition. Therefore from the previous points the sequence of mean fitnesses $(\int x p^{(\varepsilon)}_{kn+i})_{n\geq0}$ converges to $\bar{w}_i^{(\varepsilon)}$ that satisfies $\prod_{l\in\K} \bar{w}_l^{(\varepsilon)}\geq  (1-\varepsilon)^k  \prod_{l\in\K}((1-\beta_l)\eta_0)$. Finally from points $(iii)$ and $(i)$ of Proposition \ref{lem:order}, $\int x p^{(\varepsilon)}_{kn+i}(dx) \leq \int x p_{kn+i}(dx) \leq \int x p^{\bullet}_{kn+i}(dx)$ for all $n\geq 0$ and $i\in\K$. 
Hence we have the chain of inequalities
$$\bar{w}_i^{(\varepsilon)} \leq \liminf w_{kn+i} \leq \limsup  w_{kn+i} \leq \bar{w}_i^{\bullet},$$
Now we just proved that $\prod_{l\in\K} \bar{w}_l^{(\varepsilon)}\geq  (1-\varepsilon)^k \prod_{l\in\K}((1-\beta_l)\eta_0),$
while $C>0$ implies 
$ \prod_{l\in\K} \bar{w}_l^{\bullet} = \prod_{l\in\K}((1-\beta_l)\eta_0)$
by inspection of \eqref{eq:p(eta_0)}, hence 
$\prod_{l\in\K} \bar{w}_l^{(\varepsilon)} \geq (1-\varepsilon)^k \prod_{l\in\K} \bar{w}_l^{\bullet}$.
This gives, after a moment of thought,
$\bar{w}_i^{\bullet} (1- \varepsilon)^{k} \leq \bar{w}_i^{(\varepsilon)}$ for each $i \in \K$, which is enough to conclude, $\varepsilon$ being arbitrary.
\end{proof}

Proposition \ref{prop:convergence-fitness} gives the convergence of sub-sequences of mean fitnesses. We now enhance this convergence to the convergence of the sub-sequences of the fitness distributions.
We shall need the following piece of notation\footnote{beware $\bar{z}$ has nothing to do with a complex conjugate}:
\begin{equation} z_{i} := \frac{1-\beta_{[i+1]}}{\bar{w}_{i}}\text{, for } i \in \K \quad \text{ and} \quad\bar{z} := (z_{0} z_{1}\ldots z_{k-1})^{1/k},\label{eq-notation-z}\end{equation}
and we point out that the inequality \eqref{ineq:w} then translates into 
\begin{equation}
\label{ineq:w2}
\bar{z}\eta_0\leq 1.
\end{equation}
The following proposition states the convergence of the sub-sequences of fitness distributions. The proof of this result follows the strategy developed in Proposition 1 of \cite{Henard_review} and is based on the convergence of the sub-sequences of mean fitnesses established in the previous proposition.

\begin{proposition}\label{prop-cvgce-p}
\begin{itemize}
\item In the general case where $\bar{z}\eta_0\leq 1$, for each $i\in \K$, the sequence of fitness distributions $(p_{kn+i})_{n\geq0}$ converges in total variation on $[0,\xi]$ for every $\xi<\eta_0$, hence also weakly on $[0,\eta_0]$, to\footnote{the quantity $z_{[i-1]} \ldots  z_{[i-j]}$ is a shorthand for $\prod_{1 \leq \ell \leq j} z_{[i-\ell]}$, and we use the convention that products over the empty set are $1$ when $j=0$.}
\begin{equation}\label{eq-pi1}\pi_{i}(dx)=\sum_{j=0}^{k-1}\frac{z_{[i-1]} \ldots  z_{[i-j]}  }{1-(\bar{z}x)^k}x^{j} \beta_{[i-j]}q_{[i-j]}(dx)+\alpha_i\delta_{\eta_0}\end{equation} where $\alpha_i$ is the non-negative constant such that $\pi_i$ is a probability measure.
\item If furthermore $\bar{z}\eta_0<1$ then for each $i\in \K$, the sequence of fitness distributions $(p_{kn+i})_{n\geq0}$ converges in total variation on $[0,\eta_0]$, to \begin{equation*}\pi_{i}(dx)=\sum_{j=0}^{k-1}\frac{z_{[i-1]} \ldots  z_{[i-j]}  }{1-(\bar{z}x)^k}x^{j} \beta_{[i-j]}q_{[i-j]}(dx).\end{equation*}
\end{itemize}
\end{proposition}

\begin{remark}
    The inequality \eqref{ineq:w2} with condensation occurring only in the equality case is the counterpart of the floor value phenomenon for the mean limiting fitness discussed around \eqref{saturation-1d-2}.
\end{remark}

\begin{proof}[{Proof of Proposition \ref{prop-cvgce-p}, first item, local convergence in total variation on $[0,\eta_0)$, and $\bar{z}\eta_0\leq 1$}]
Let $\xi <\eta_0$.
We wish to compare:
\begin{align*}
p_{kn+i}(dx) 
& = \sum_{j=0}^{kn+i-1}
\underbrace{\left(\prod_{\ell=0}^{j-1}\frac{1-\beta_{[i-\ell]}}{w_{kn+i-\ell-1}}\right) x^j\beta_{[i-j]}q_{[i-j]}(dx)}_{a_{kn+i,j}(dx)}
+ \underbrace{\left(\prod_{\ell=0}^{kn+i-1}\frac{1-\beta_{[i-\ell]}}{w_{kn+i-\ell-1}}\right) x^{kn+i} p_0(dx)}_{b_{kn+i}(dx)}\\
\end{align*}
and 
\begin{align*}
\pi^{\circ}_{i}(dx) 
& = \sum_{j=0}^{k-1}\frac{z_{[i-1]} \ldots  z_{[i-j]}  }{1-(\bar{z}x)^k}x^{j} \beta_{[i-j]}q_{[i-j]}(dx)  = \sum_{j=0}^{\infty} \underbrace{\left(\prod_{\ell=0}^{j-1}\frac{1-\beta_{[i-\ell]}}{\bar{w}_{[i-\ell-1]}}\right) x^j\beta_{[i-j]}q_{[i-j]}(dx)}_{\bar{a}_{i,j}(dx)} 
\end{align*}
in the sense of the total variation distance on $[0,\xi]$.
We point out that Fatou Lemma ensures that $\pi^{\circ}_i(dx)$ is a sub-probability measure. 
First, 
\begin{align*}
 \|a_{kn+i,j}(dx) -  \bar{a}_{i,j}(dx) \|_{TV, [0,\xi]}  
& = \left\lvert \prod_{\ell=0}^{j-1}\frac{1-\beta_{[i-\ell]}}{w_{kn+i-\ell-1}}  
- \prod_{\ell=0}^{j-1}\frac{1-\beta_{[i-\ell]}}{\bar{w}_{[i-\ell-1]}}  \right\lvert \int_{0}^{\xi} x^j\beta_{[i-j]}q_{[i-j]}(dx)\\
& = \left\lvert 
\prod_{\ell=0}^{j-1}\frac{\bar{w}_{i-\ell-1}}{w_{kn+i-\ell-1}}  
-1  \right\lvert  \bar{a}_{i,j}([0,\xi])  
\end{align*}
hence, for any $r_0$:  
\begin{align}
\label{1a}
\sum_{j=0}^{kr_0-1} \|a_{kn+i,j}(dx) -  \bar{a}_{i,j}(dx) \|_{TV, [0,\xi]}
\leq   \sum_{j=0}^{kr_0-1} \left\lvert 
\prod_{\ell=0}^{j-1}\frac{\bar{w}_{i-\ell-1}}{w_{kn+i-\ell-1}}  
-1  \right\lvert  \bar{a}_{i,j}([0,\xi])  
\leq \max_{0 \leq j < kr_0}\left\lvert 
\prod_{\ell=0}^{j-1}\frac{\bar{w}_{i-\ell-1}}{w_{kn+i-\ell-1}}  
-1  \right\lvert   
\end{align}
using $\sum_{j=0}^{\infty } \bar{a}_{i,j}([0,\xi]) \leq 1$ for the last inequality, and $r_0$ being fixed, the last term goes to $0$ as $n \to \infty$, using the convergence $w_{kn+i-\ell-1} \to \bar{w}_{i-\ell-1}$ established in Proposition \ref{prop:convergence-fitness}.

Second, we notice $\bar{z} \xi<1$ by assumption and write:
\begin{align}
\sum_{j \geq  kr_0} \bar{a}_{i,j}([0,\xi]) 
& = \int_{0}^{\xi}  \sum_{j \geq  kr_{0}} 
\left( \prod_{\ell=0}^{j-1}\frac{1-\beta_{[i-\ell]}}{\bar{w}_{[i-\ell-1]}} \right) x^j \beta_{[i-j]} q_{[i-j]}(dx)  \nonumber \\
& = (\bar{z} \xi)^{kr_0}  \int_{0}^{\xi}  \sum_{j \geq  0} 
\left( \prod_{\ell=0}^{j-1}\frac{1-\beta_{[i-\ell]}}{\bar{w}_{[i-\ell-1]}} \right) x^j \beta_{[i-j]} q_{[i-j]}(dx)  \nonumber \\
&  \leq    (\bar{z} \xi)^{kr_0}  \sum_{j \geq  0} \bar{a}_{i,j}([0,\xi]) \nonumber \\
& \leq    (\bar{z} \xi)^{kr_0}.  
 \label{2a} 
\end{align}

Third, choosing $\delta>0$ small enough such that
$\bar{z} \xi(1+\delta) <1$, and then $n_0$ large enough such that
$(1+\delta) \prod_{\ell=0}^{k-1} w_{kn+i-\ell-1} \geq  \prod_{i \in \K}  \bar{w}_{i}$ for $n \geq n_0$, then for $n$ such that $n \geq n_0 +r_0$, we find:
\begin{align}
\sum_{j =k r_0}^{kn+i-1} a_{kn+i,j}([0,\xi])  
 & = \int_{0}^{\xi}  \sum_{j = k r_{0}} ^{kn+i-1 } \left(\prod_{\ell=0}^{j-1}\frac{1-\beta_{[i-\ell]}}{w_{kn+i-\ell-1}}\right) x^j\beta_{[i-j]}q_{[i-j]}(dx) \nonumber\\
 & \leq  \left(\bar{z} \xi (1+\delta)\right)^{kr_0}\int_{0}^{\xi}  \sum_{j = 0} ^{k(n-r_0)+i-1 } \left(\prod_{\ell=0}^{j-1}\frac{1-\beta_{[i-\ell]}}{w_{k(n-r_0)+i-\ell-1}}\right) x^j\beta_{[i-j]}q_{[i-j]}(dx) \nonumber \\ 
& = \left(\bar{z} \xi (1+\delta)\right)^{kr_0}   \sum_{j =0}^{k(n-r_0)+i-1} a_{k(n-r_0)+i,j}([0,\xi])  
  \nonumber \\ 
  & \leq  \left(\bar{z} \xi (1+\delta)\right)^{kr_0} .\label{3a}
\end{align}

Fourth, the term implying $b_{kn+i}(dx)$ is dealt with as the second term above: 
\begin{align}
b_{kn+i}([0,\xi])
&  = \int_{0}^{\xi} \sum_{j=kr_{0}}^{kn+i-1 } \left(\prod_{\ell=0}^{j-1}\frac{1-\beta_{[i-\ell]}}{w_{kn+i-\ell-1}}\right) x^j\beta_{[i-j]}q_{[i-j]}(dx) \nonumber \\
& \leq (\bar{z} \xi (1+\delta))^{kr_0} \int_{0}^{\xi} \sum_{j=0}^{k(n-r_0)+i-1 } \left(\prod_{\ell=0}^{j-1}\frac{1-\beta_{[i-\ell]}}{w_{k(n-r_0)+i-\ell-1}}\right) x^j\beta_{[i-j]}q_{[i-j]}(dx) \nonumber \\
& = \left(\bar{z} \xi (1+\delta)\right)^{kr_0}  
   b_{k(n-r_0)+i}([0,\xi])  
  \nonumber \\ 
&\leq  (\bar{z} \xi (1+\delta))^{kr_0}.  \label{4a}
\end{align}
The sum of the four terms in \eqref{1a},\eqref{2a}, \eqref{3a} and \eqref{4a} now gives an upper bound for $\|p_n-\pi\|_{TV, [0,\xi]}$.
Now, we choose the parameters as follows: we choose $r_0$ large enough such that \eqref{3a}\eqref{4a} (hence \eqref{2a}) are small and then $n \geq n_0+r_0$ large enough so as to make \eqref{1a} small.
\end{proof}

\begin{proof}[{Proof of Proposition \ref{prop-cvgce-p}, first item, weak convergence on {$[0,\eta_0]$}, and $\bar{z}\eta_0\leq 1$}]
Fix $i \in \K$. The set of probability measures on $[0,\eta_0]$ is compact for the topology of weak convergence; now take any subsequence of $(p_{kn+i}(dx))_{n \geq 0}$; it admits a converging subsubsequence, which should agree with the limit in total variation on $[0,\eta_0)$, that is with  $\pi^{\circ}_i$, and be supported on $[0,\eta_0]$: the only possibility left is then given  by $\pi_i$ defined in \eqref{eq-pi1} on $[0,\eta_0]$. The set of accumulation points of $(p_{kn+i}(dx))_{n \geq 0}$ therefore consists in the singleton given by $\{\pi_i\}$: in other words, the sequence $(p_{kn+i})_{n \geq 0}$ has weak limit $\pi_i$ on $[0,\eta_0]$.
\end{proof}

\begin{proof}[{Proof of Proposition \ref{prop-cvgce-p}, second item, convergence in total variation on {$[0,\eta_0]$}, and $\bar{z}\eta_0< 1$}]

The proof is similar, even easier, and consists in replacing $\xi$ by the larger quantity $\eta_0$ in those bounds involving $\xi$; precisely, \eqref{1a} is unaffected, then choosing 
$\delta$ such that $\bar{z}\eta_0(1+\delta) <1$ (and $n_0$ as before), we get the following substitutes for \eqref{2a} \eqref{3a} \eqref{4a}:
\begin{align*}
& \sum_{j \geq  kr_0} \bar{a}_{i,j}([0,\xi]) 
 \leq    (\bar{z} \eta_0)^{kr_0} \\
& \sum_{j \geq  kr_0} a_{kn+i,j}([0,\xi]) 
 \leq    (\bar{z} \eta_0(1+\delta))^{kr_0} \\
&  b_{kn+i,j}([0,\xi]) 
 \leq    (\bar{z} \eta_0 (1+\delta))^{kr_0}.
\end{align*}

\end{proof}

\subsection{Characterization of the limit and proof of Theorem \ref{main2}}
Having proved in the previous section that the sequence $((p_{kn+i}(dx))_{i \in \K})_{n}$ converges, we now characterize the limiting distribution $(\pi_{i}(dx))_{i \in K}$ by identifying its mean vector $(\bar{w}_i)_{i \in \K}$
as the solution of a fixed point equation; along the way, we find the criterion for condensation in terms of a function of the real variable encapsulating the non-linearity of the problem.

Recall the definition of $(z_i)_{i \in \K}$ and $\bar{z}$ in \eqref{eq-notation-z} and define the column vector $V=(v_0,...,v_{k-1})\in(\mathbb{R}_+^{\star})^k$ by 
\begin{equation}
    v_{i}:= \frac{z_{k-1} z_{k-2} \ldots z_{i}}{(\bar{z})^{k-i}} ,  \quad \text{ so that } v_0 =1
    \label{def-v}
\end{equation}
and let $\bf 1$ be the column vector with coordinates $(1, \ldots, 1)$. The following proposition shows that the vector $V$ solves a linear system.

\begin{proposition}\label{prop-equation}
It holds:
$$ A(\bar{z}) \, V+ \alpha \mathbf{1} = V  $$ for a non-negative constant $\alpha \geq 0$.
Furthermore, in case $\bar{z} < 1/\eta_0$, it holds $\alpha = 0$.
\end{proposition}

\begin{proof}
Recall  we defined the measure $\pi^{\circ}_i(dx)$, that one may think of as $\pi_i(dx)$  possibly deprived from its mass at $\eta_0$, by
\begin{align*}
 \pi^{\circ}_i(dx) &:= \sum_{j=0}^{k-1} \Big(\frac{z_{[i-1]} \ldots  z_{[i-j]}}{1- (\bar{z} x)^k}\Big) x^j \beta_{[i-j]} q_{[i-j]}(dx), \\
&= \sum_{j=0}^{k-1} \Big(\frac{z_{[i-1]} \ldots  z_{[i-j]}}{\bar{z}^j}\Big)  \frac{(\bar{z} x)^j \beta_{[i-j]} q_{[i-j]}(dx)}{1- (\bar{z} x)^k} \nonumber
\end{align*}
and observe, using the definition \eqref{eq:mu-ij} of $\mu_i^{j}(z)$, that one has:
\begin{align*}
 \int \pi^{\circ}_i(dx) 
&= \sum_{j=0}^{k-1} \Big(\frac{z_{[i-1]} \ldots  z_{[i-j]}}{\bar{z}^j}\Big)  \mu_{[i-j]}^j(\bar{z})\nonumber
\end{align*}
We already know from Proposition \ref{prop-cvgce-p} that $\pi^{\circ}_i(dx)$ is a sub-probability measure in general.
Now, in case  $\pi^{\circ}_i(dx)$ is a probability measure,
multiplying the last display by $v_{i}:= \frac{z_{k-1} \ldots z_{i}}{(\bar{z})^{k-i}}$  gives an equation involving the quantities $(v_{i})_{i \in \K}$ only (use the key fact that $v_0=1$): 
$$\sum_{j=0}^{k-1}  \mu_{[i-j]}^{j}(\bar{z})  v_{[i-j]}  = v_{i} $$
which, upon a change of the indices gives
$\sum_{\ell=0}^{k-1}  \mu_{\ell}^{[i-\ell]}(\bar{z})  v_{\ell}  = v_{i} $, that is  $A(\bar{z}) \; V = V$. In the general case, $\pi^{\circ}_i(dx)$ is a sub-probability measure  only which, following the same arguments as before, translates into the (coordinate-wise) inequality:
$A(\bar{z}) \; V \leq V$.
In this case, our way to an equality is to observe the following key relation:
$$\int \Big(1-\frac{x}{\eta_0}\Big) \pi^{\circ}_i(dx) = \int \Big(1-\frac{x}{\eta_0}\Big) \pi_i(dx) = 1-\frac{w_i}{\eta_0},$$
where the first equality follows from the fact that $\pi^{\circ}_i(dx)$ and $\pi_i(dx)$ agree on $[0,1/\eta_0)$. Again, we translate this equation in terms of the matrix $A(\bar{z})$. First we have:
\begin{equation}\label{eq-(1-x)nu}\int \Big(1-\frac{x}{\eta_0}\Big) \pi^{\circ}_i(dx) = \sum_{j=0}^{k-1} \frac{z_{[i-1]} \ldots  z_{[i-j]}}{(\bar{z})^j}  \int \Big(1-\frac{x}{\eta_0}\Big) \frac{(\bar{z} x)^j \beta_{[i-j]} q_{[i-j]}(dx)}{1- (\bar{z} x)^k},\end{equation}
and the relation  
$\mu^k_\ell(z)-\mu^0_\ell (z) = -\beta_{\ell}$
then gives
(pay attention that $j+1$ is under bracket in the last term):
$$
\int_{0}^{1} x \frac{(\bar{z}x)^j \beta_{[i-j]} q_{[i-j]}(dx)}{1-(\bar{z}x)^k}  
 = \frac{1}{\bar{z}} \mu_{[i-j]}^{j+1}(\bar{z})   
 = \frac{1}{\bar{z}} \left(\mu_{[i-j]}^{[j+1]}(\bar{z}) - \Ind{j+1=k}   \beta_{[i+1]}\right)
$$
so that the RHS of \eqref{eq-(1-x)nu} writes
$$\sum_{j=0}^{k-1} \frac{ z_{[i-1]} \ldots z_{[i-j]}}{(\bar{z})^{j}} \left(\mu_{[i-j]}^{j}(\bar{z})-  \frac{1}{\eta_0 \bar{z}} \left( \mu_{[i-j]}^{[j+1]}(\bar{z}) - \Ind{j+1=k}   \beta_{[i+1]} \right) \right) $$
while the LHS is in our notation:
$$1-\frac{w_i}{\eta_0} = 1 - \frac{1-\beta_{[i+1]}}{\eta_0 z_i}.$$
Multiplying these two equal expressions by $v_{i}$ defined in \eqref{def-v} gives:
$$\sum_{j=0}^{k-1} \left(\mu_{[i-j]}^{j}(\bar{z})- \frac{1}{\eta_0 \bar{z}} \left( \mu_{[i-j]}^{[j+1]}(\bar{z}) - \Ind{j=k-1}   \beta_{[i+1]} \right)\right) v_{[i-j]}  = v_i - \frac{1-\beta_{[i+1]}}{\eta_0 \bar{z}} v_{[i+1]} $$
Reindexing and observing the terms in factor of $\beta_{[i+1]}$ are the same on the LHS and the RHS, so that 
$$\sum_{\ell=0}^{k-1}  \mu_{\ell}^{[i-\ell]}(\bar{z})  v_{\ell} - v_{i}   =
 \frac{1}{\eta_0 \bar{z}} \left( \left(\sum_{\ell=0}^{k-1}  \mu_{\ell}^{[i+1-\ell]}(\bar{z})\right)  v_{\ell} -  v_{[i+1]} \right) $$
which is the equation:
$$(A(\bar{z}) V - V)_i =\frac{1}{\eta_0 \bar{z}} (A(\bar{z}) V - V)_{[i+1]} .$$
Now, in case $\eta_0 \bar{z}<1$, iterating this relation $k$ times, we recover the equation $A(\bar{z}) V =V$ \footnote{but we feel there is some pedagogical value in keeping the first proof of this relation displayed above}, while in case  $\eta_0 \bar{z}=1$, we obtain:
$A(\bar{z}) V +\alpha \mathbf{1}=V$ for a certain constant $\alpha$, that we already know to be non-negative, since $A(\bar{z}) V  \leq V$ coordinate-wise.
\end{proof}

\begin{lemma}
\label{prop-A}
The Perron eigenvalue $\rho(A(z))$ of the non-negative matrix $A(z)$ satisfies the following:
\begin{enumerate}
\item $\rho(A(0)) = \max\{ \beta_i, i \in \K\}<1$
\item $z \mapsto \rho(A(z))$ is continuous and increasing on $I_{\mu}$.
\end{enumerate}
\end{lemma}

\begin{proof}[Proof of Lemma 2]
By definition, the matrix $A(0)$ is diagonal with non-negative entries $\beta_0, \ldots, \beta_{k-1}$ on the diagonal, hence its Perron-Frobenius eigenvalue is $\max\{\beta_i, i \in \K\}$, as claimed in point 1.
Also, the fact that $z \mapsto \rho(A(z))$ is increasing follows from the fact the maps $z \mapsto A_{i,j}(z)$ are  increasing for every $i,j \in \K$ and Proposition \ref{prop-Perron} (iv) in Appendix \ref{sec-appendix}, while the continuity of $z \mapsto \rho(A(z))$ is a consequence of the continuity of the maps $z \mapsto A_{i,j}(z)$ and Proposition \ref{prop-Perron} (iii).
\end{proof}

We are now in position to prove Proposition \ref{prop:2cases}.
We start with a clarification:
\begin{lemma}
\label{1eta0-def}
$(1/\eta_0 \notin I_{\mu})$ implies $\rho(A((1/\eta_0)-))=+\infty$.
\end{lemma}

\begin{proof}
In case $1/\eta_0 \notin I_{\mu}$, necessarily $\eta_q= \eta_0$ and $\rho_{\ell}((1/\eta_q)-)$ is infinite for some $\ell \in \K$ by definition, which also implies $A_{i,\ell}((1/\eta_q)-)= \mu_{\ell}^{[i-\ell]}((1/\eta_q)-) =+ \infty$ for each $i \in \K$. Now, the row sums bound given in Proposition \ref{prop-Perron} point (ii) gives: 
$$+\infty = \min_{i} A_{i,\ell}((1/\eta_q)-) \leq \min_{i} \sum_{j \in \K} A_{ij}((1/\eta_q)-) \leq \rho(A((1/\eta_q)-)).$$
\end{proof}

\begin{proof}[Proof of Proposition \ref{prop:2cases}] 
Lemma \ref{prop-A} ensures that the coefficients of the matrix $A(z)$ are continuous, increasing functions of $z \in I_{\mu}$ and that $\rho(A(0))<1$, hence the following alternative holds: 
\begin{enumerate}
    \item If $\rho(A(1/\eta_0))\geq 1$, there exists a unique $z_c\in ]0,1/\eta_0]$ (0 is excluded because $\rho(A(0))<1$ by point 1 of Lemma \ref{prop-A}) such that $\rho(A(z_c))=1$, in which case the Perron theorem for matrices with positive coefficients ensures there exists a unique vector $X \in \, (\mathbb R_+^\star)^k$ such that $A(z_c)X=X$ and $X_1=1$; conversely, a vector $X$ with positive coefficients solving the equation $A(z) X = X$ for some $z>0$ is necessarily the Perron eigenvector hence $1$ is the Perron eigenvalue of $A(z)$ which implies $z=z_c$. 
    \item If $\rho(A(1/\eta_0))<1$, then by Lemma \ref{1eta0-def}, $1/\eta_0 \in I_{\mu}$ so the matrix $A(1/\eta_0)$, that we simply write $A$ here,  is well defined and the eigenvalues of $A$ (in the complex plane $\mathbb C$) belong to $\{z \in \mathbb C  : |z| \leq  \rho(1/\eta_0)\}$ hence the eigenvalues of the matrix $I-A$ belong to $\{z \in \mathbb C  : |1-z| \leq  \rho(1/\eta_0)\}$, and this set does not contain $0$: this means the matrix $I-A$ is invertible, 
    and one can furthermore write its inverse as
    $(I-A)^{-1} =\sum_{n\geq 0} A^n$. Hence there exists a unique $Y\in \R^k$ such that $(I-A) Y=\mathbf{1}$ and $X= Y/Y_0$ satisfies 
    $(I-A)X=(1/Y_0) \mathbf{1}$ and $X_0 =1$. Last
    $X= (\sum_{n \geq 0} A^n) (1/Y_0) \mathbf{1}$ so that $\alpha=1/Y_0$ and the coordinates of $X$ indeed have the same sign, hence are all positive since $X_0 = 1$, which in turn implies $\alpha>0$. 
\end{enumerate} 
\end{proof}

\begin{proof}[Proof of Theorem \ref{main2}]
We have proved that for each $i \in \K$, the sequence $(p_{kn+i}(dx))_{n \geq 0}$ converges in the appropriate sense towards a probability measure
$\pi_i(dx)$ parameterized by coefficients $(z_i)_{i \in \K}$ (Proposition \ref{prop-cvgce-p}). Next we have proved that these coefficients $(z_i)_{i \in \K}$ are solution of a fixed point equation (Proposition \ref{prop-equation}) that agrees with \eqref{U-theta}. Finally we have proved that Equation \eqref{U-theta} admits a unique solution (Proposition \ref{prop:2cases}), so that $\bar{z}$ and $V$ respectively defined in \eqref{eq-notation-z} and \eqref{def-v} indeed agree with $z_c$ and $U$ in the statement of Theorem \ref{main2}.
\end{proof}

\begin{proof}[Proof of Proposition \ref{closed-formula}] 
Summing the coordinates of the equality \eqref{U-theta}, $A(z_c)U+ \alpha \mathbf{1} = U$, as in \eqref{sum:columns}
gives:
$$\alpha=\frac{1}{k}\sum_{i=0}^{k-1}U_i\left(1- \rho_i(z_c) \right).$$  
After replacing $\alpha$ by this formula, Equation \eqref{U-theta} writes down:
$BU=0$, for $B= B(z_c)$, and we establish in a first part of the proof that the vector $(N_i)_{i\in\K}$ also solves $BN=0$.

Let $M_{ij}$, $i,j \in \K$ be the $(i,j)$-minor of the matrix $B$ (the determinant of the matrix $B$ with line $i$ and column $j$ both removed). We make the claim that these minors satisfy \begin{equation}\label{minors-rank-1}M_{ij}=(-1)^{i+j}M_{jj}\end{equation} for all $i,j\in\K$. 
To prove Equation \eqref{minors-rank-1}, we observe that each column of $B$ has a null sum, which follows from \eqref{sum:columns}, then perform standard manipulations on the matrix $B$:
\begin{align*}
M_{ij}
&=\det\Big((B_{lm})_{\substack{0\leq l,m \leq k-1\\l\neq i, m\neq j}}\Big)=\det\Big((B_{lm}+\sum_{r\neq i, r\neq j} B_{rm}\mathbf{1}_{\ell=j})_{\substack{0\leq l,m \leq k-1\\l\neq i, m\neq j}} \Big)\\
&=\det\Big((B_{lm}\mathbf{1}_{l\neq j}-B_{im}\mathbf{1}_{l=j})_{\substack{0\leq l,m \leq k-1  \\ l\neq i, m \neq j }}\Big)=(-1)^{i+j} \det((B_{lm})_{\substack{0\leq l,m \leq k-1  \\ l\neq j, m \neq j }})=(-1)^{i+j}M_{jj},\end{align*}
adding up (second equality) all the other lines to the line indexed by $j$, which gives (third equality) the negative of the $i$-th line of the original matrix $B$ since each column of that matrix sums to $0$, then observing that the resulting matrix is the one with the $j$-th line erased, up to permutations of the lines that result in the factor $(-1)^{i+j}$ (fourth equality).
Using the claim, we finally observe that computing the $i$-th coordinate $(BN)_i$ amounts to expanding $\det(B)$ according to the $i$-th line:
$$(BN)_i=\sum_{j\in\K} B_{ij}N_{j}=\sum_{j\in\K} B_{ij}M_{jj}=\sum_{j\in\K} (-1)^{i+j}B_{ij}M_{ij}=\det(B)=0,$$
where the last equality follows again from the fact that the sum of the columns of $B$ is null. Equation \eqref{minors-rank-1} is now proved. We notice for future use that it implies the rank of the adjugate matrix adj$(B)$ of $B$ (see Proposition \ref{adjugate} for its definition and a few properties) is $\leq 1$.

Now for the second part of the proof. 
The equation $BN=0$ implies that $A N + \alpha' 1 = N$ holds for some constant $\alpha'$ and then:
\begin{enumerate}
\item In case $\rho(A(z_c))<1$, Proposition \ref{prop:2cases} point 2 ensures that $(N,\alpha')$ and $(U,\alpha)$ are proportional, and that $B$ has rank $k-1$; assume $(N,\alpha')$ is identically null; in such a case, the adjugate matrix adj$(B)$ being of rank $\leq 1$ with a null diagonal $N$ is then null, which in turn implies $B$ is of rank $\leq k-2$ (Proposition \ref{adjugate}), a contradiction.
\item In case $\rho(A(z_c))=1$, 
if $\alpha'\neq 0$, one can assume $\alpha'>0$ wlog, and then, choosing $\gamma>0$ large enough, the vector $X= N + \gamma U$ has positive coordinates (because $U$ has by Perron Theorem) and satisfies $A(z_c)X + \alpha' \mathbf{1} = X$, in particular $A(z_c)X>X$, which cannot hold by Collatz-Wielandt formula (Proposition \ref{prop-Perron} $(i)$); therefore $\alpha'=0$, and $N$ and $U$ are proportional by Perron Theorem; as above, $N$ being identically null would imply $B$ has rank $\leq k-2$, hence we could find two linearly independent vector in the kernel of $B$, that is, by the above reasoning about the nullity of $\alpha'$, two linearly independent solutions $X \in \R^k$ to the equation $AX=X$: this cannot hold by Perron Theorem again.
\end{enumerate}
In both cases, $N$ is non identically null and proportional to $U$ with positive coordinates, hence $N_0 \neq 0$ and $U_i= N_i/N_0, i \in \K$ follows using $U_0 =1$.

\end{proof}

\section{Some additional results}
\label{sec-additional}
\subsection{Alternative form of the criterion}\label{sec-alternative}

We aim in this section at an alternative formulation of our criterion. Recall $\eta_q = \max\{\eta_{q_i}, i \in \K \}$, and introduce the following determinant:
$$\Psi : I_{\mu} \to \R, z \mapsto \Psi(z)=\det(\mathbf{I}_k-A(z)),$$
where $\mathbf{I}_k$ stands for the identity matrix of dimension $k$. Our goal is to prove Proposition \ref{equivalence2conditions},  that states that $1-\rho(A(z))$ and $\Psi(z)$ have the same sign for $z \in I_{\mu}$.
 .
We show that the statement is a corollary of the following matrix-theoretic result proved in the appendix, one of our main efforts in this paper.
\begin{proposition}
\label{at-most-one}
For $ z \in I_{\mu}$, $A(z)$ has at most one real eigenvalue that is $\geq 1$. (In case it has one, this eigenvalue is then necessarily its Perron eigenvalue $\rho(A(z))$ with multiplicity $1$).
\end{proposition}

\begin{proof}[Proof of Proposition \ref{equivalence2conditions}]
Setting $\{\lambda_i(z), i \in \K\}$ to be the set of eigenvalues of $A(z)$ counted with multiplicity, we have:
$$\Psi(z) = \prod_{i \in \K} (1 -\lambda_i(z)).$$
In this product: 
\begin{enumerate}
    \item the real eigenvalues <1 contribute to a positive term.
    \item the complex eigenvalues that are not real contribute to a positive term, since each such eigenvalue can be paired to its conjugate value.
\end{enumerate}
Now, $\rho(A(z))$ is the only real eigenvalue that may be $\geq 1$ by Proposition \ref{at-most-one}, also, it has multiplicity $1$, hence $\Psi(z)$ is of the sign of $1- \rho(A(z))$.
\end{proof}

$\Psi$ defines a continuous map on $I_{\mu}$ that, on the basis of Proposition \ref{equivalence2conditions} and the properties of $\rho \circ A$ recalled in Lemma \ref{prop-A}, satisfies the following: as $z$ increases from $0$, $\Psi$ is  first positive and, then depending on  the sign of $\rho(A((1/\eta_q)-))$, it may cancel at a single point (since $\rho$ is increasing), and then becomes negative. 
(Let us point though that contrarily to $\rho \circ A$, the map $\Psi$ has no a-priori monotonicity properties).\\

In the case with $k=2$ environments, we find an expression that is an even simpler expression than $\Psi(z)$ and shares the same sign, namely:
$$\Gamma_2(z) := \frac{1- \int  \frac{\beta_0 q_0(dx)}{1- zx} }{1- \int  \frac{\beta_0 q_0(dx)}{1+ zx}  }+ \frac{1- \int  \frac{\beta_1 q_1(dx)}{1- zx} }{1- \int  \frac{\beta_1 q_1(dx)}{1+ zx }},$$
where we stress that the two denominators are positive (since $\beta_0,\beta_1$ are both $<1$). The sign of this expression gives an alternative way to express the criterion for condensation.
\begin{proposition}\label{prop-Gamma2}
Let $k=2$. For $z \in I_{\mu}$, the two expressions $\Gamma_2(z)$ and  $1-\rho(A(z))$ have the same sign, in the sense they are both negative, null or positive.
\end{proposition}

Notice that $\Gamma_2(z)$ is sum of two terms, each depending only on a single mutation environment; one can not hope to find a quantity with the same form in case $k\geq 3$ since the model is then no more invariant by permutation.

\begin{proof} 
First, observe that, for $\varepsilon \in \{-1,1\}$, $i \in \K$ and $z\leq 1/\eta_{q_i}$,
$$\mu_i^{0}(z) + \varepsilon \mu_i^{1}(z) 
= 
\int  \frac{ \beta_i q_i(dx)}{1- (zx)^2} 
+ \varepsilon \int  \frac{ zx \beta_i q_i(dx)}{1- (zx)^2}
= \int  \frac{ \beta_i q_i(dx)}{1 - \varepsilon zx}. 
$$ 
Therefore $\Gamma_2(z)$ has the same sign than the following expression:
\begin{align*} &  \Big(1- \int  \frac{\beta_0 q_0(dx)}{1- zx}\Big)\Big(1- \int  \frac{\beta_1 q_1(dx)}{1+ zx }\Big)+\Big(1- \int  \frac{\beta_1 q_1(dx)}{1- zx} \Big)\Big(1- \int  \frac{\beta_0 q_0(dx)}{1+ zx} \Big)\\
 & = \Big(1-(\mu_0^{0}(z)+\mu_0^{1}(z))\Big) \Big(1-(\mu_1^{0}(z)-\mu_1^{1}(z))\Big) \\
&+\Big(1-(\mu_0^{0}(z)-\mu_0^{1}(z))\Big) \Big(1-(\mu_1^{0}(z)+\mu_1^{1}(z))\Big)\\
 & = 2 \Big((1-\mu_0^{0}(z))(1-\mu_1^{0}(z))-\mu_0^{1}(z)\mu_1^{1}(z)\Big)\\
& = 2 \det(I-A(z)) = 2 \Psi(A(z))
 \end{align*}
and the claim now follows from Proposition \ref{equivalence2conditions}.
\end{proof}

\subsection{Generating function of the weights}

In this section, we give some closed formulas for the generating functions of cumulative products of the mean fitnesses. These cumulative products, called weights, are instrumental in the proof by Kingman. They could also be used here  to prove the convergence of the fitness along subsequences, and even to get bounds on the speed of convergence, in the case $\bar{z}<1/\eta_0$; alas, in case
$\bar{z}=1/\eta_0$, the ingenious argument of Kingman based on a certain completely monotonous sequence does not seem to extend to our setting of periodic mutation environments.
We still present these formulas for two reasons: they have interest by themselves, and they help to justify why the determinant $\Psi(z)$ is a natural quantity for the problem.

We now set up some notation. First, let $\bar{\beta}$ be the unique non-negative real number satisfying $(1-\bar{\beta})^k =\prod_{i=0}^{k-1} (1-\beta_i)$, and, for $i,j \in \K$, set:
$$ c_{i,j}=\begin{cases}\frac{\prod_{q=0}^{[i-j]-1}(1-\beta_{[i-q]})}{(1-\bar{\beta})^{[i-j]}} & \text{if $ i\neq j$}\\1 & \text{if $i=j$.}\end{cases}
$$
The numbers $c_{i,j}$ and  $c_{j,i}$ are inverse of each other: $c_{i,j} c_{j,i}=1$ for $0 \leq i,j <k$. We recall the definition \eqref{eq:mu-ij} of (a variant of) the moment generating generating function $\mu_{i}^{j}(z)$ of the sub-probability  measure $\beta_i q_i(dx)$ and also introduce its $n$-th moment $\mu_{n,i}$:
$$\mu_{i}^{j}
(z) = \sum_{n \geq 0 : [n]=j} z^{n} \mu_{n,i}= \int_{0}^{1} \frac{(zx)^{j}}{1-(zx)^k} \beta_i q_i(dx),  \qquad   i,j \in \K. $$
Now let $W_n = \prod_{0 \leq k < n} w_k$ be the product of the mean fitnesses, and introduce the following generating functions, again distinguished by their rest in the division by $k$, and set:
$$  w^{(i)}(z)  = \sum_{n \geq 1 : [n]=i}\frac{W_{n}}{(1-\bar{\beta})^{n}} z^{n}, \qquad  i \in K.
$$
Last, we shall need $m_n=\int x^n p_0(dx)$ and: 
$$m^{(i)}(z) = \sum_{n\geq 0 : [n]=i} z^{n} m_{n} =  \int_{0}^{1} \frac{(zx)^i}{1-(zx)^k} p_0(dx), \qquad   i \in K.$$
\begin{proposition}
For $z$ a complex number such that $|z| < 1/\eta_0$:
\begin{equation}
\label{eq:k-key}
 w^{(i)}(z)  = \sum_{j=0}^{k-1} c_{i,j} 
 \mu_{j}^{[i-j]}(z) w^{(j)}(z)
   + c_{i,0} m^{(i)}(z).
\end{equation}
\end{proposition}

\begin{proof}
We start with a general inhomogeneous environment described by a sequence of sub-probability measures  $(\beta_1 q_1, \beta_2 q_2, \ldots )$: 
\begin{equation*} p_{n+1}(dx)  = \beta_{n+1} q_{n+1}(x) + (1-\beta_{n+1}) \frac{x p_{n}(dx)}{w_{n}}, \qquad n \geq 0.\end{equation*}
In this setting, an induction over $n$ gives the following representation of $p_{n+1}$ as a function of the initial distribution and the mutations environments $(\beta_1 q_1, \beta_2 q_2, \ldots )$: 

  $$p_{n+1}(dx) =\sum_{\ell=0}^{n} \Big(\prod_{n+2-\ell \leq  j  \leq n+1} \frac{1- \beta_j}{w_{j-1}}\Big)    x^\ell  \beta_{n+1-\ell} q_{n+1-\ell}(dx) + \Big(\prod_{1 \leq j  \leq n+1}  \frac{1- \beta_j}{w_{j-1}}\Big) x^{n+1} p_0(dx).
  $$

Multiplying by $W_{n+1} = \prod_{0 \leq j \leq n} w_j$ and integrating, we get:
$$W_{n+1} = \sum_{\ell=0}^{n} \Big(\prod_{j=n+2-\ell}^{n+1} (1- \beta_j) \Big) W_{n+1-\ell}  \mu_{\ell,n+1-\ell} 
+  \prod_{j=1}^{n+1} (1- \beta_j) \; m_{n+1}.$$

In particular, 
\begin{equation*}
W_{kn+i} =
\sum_{l=0}^{kn+i-1} \Big(\prod_{j=kn+i+1-l}^{kn+i} (1- \beta_j) \Big)W_{kn+i-l} \mu_{l,kn+i-l} 
+  \prod_{j=1}^{kn+i} (1- \beta_j) m_{kn+i}.
\end{equation*}
We now apply the formula in our periodic setting where $\beta_n  q_n= \beta_{[n]}q_{[n]}$. Using the notation $V_n=W_n/(1-\bar{\beta})^n$ with $\bar{\beta}$ defined by $(1-\bar{\beta})^k =\prod_{j=0}^{k-1} (1-\beta_j)$, and  distinguishing according to the remainder in the division by $k$, we get for any $i\in\K$ (void products are to be treated as $1$): 
\begin{align*}
V_{kn+i} & = \sum_{j=0}^{k-1}
\frac{\prod_{q=0}^{j-1} (1- \beta_{[i-q]})}{(1-\bar{\beta}^{j})} \sum_{0 \leq k\ell+j \leq kn+i-1} V_{kn+i- (k\ell+j)} \mu_{k\ell+j,[i-j]}  + \frac{\prod_{q=1}^{i} (1- \beta_{q})}{(1-\bar{\beta})^{i}}m_{kn+i}
\end{align*}
which entails the following for the generating functions $w^{(i)}$:
$$
w^{(i)}(z)  =\sum_{j=0}^{k-1}
\frac{\prod_{q=0}^{j-1} (1- \beta_{[i-q]})}{(1-\bar{\beta})^{j}} 
w^{([i-j])}(z)
 \mu_{[i-j]}^{(j)}(z)  + \frac{\prod_{q=1}^{i} (1- \beta_{q})}{(1-\bar{\beta})^{i}}m^{(i)}(z)
$$
as claimed in our proposition in our newly set notation.
\end{proof}

We now set matrices notations:
$$A= (\mu_j^{[i-j]}(z))_{0 \leq i,j <k}, \qquad  C = (c_{i,j})_{0 \leq i,j <k}, \qquad M = (m^{(i)})_{0 \leq i < k}, \qquad  W=(w^{(i)}(z))_{0\leq i<k}, $$
also, we let $C_0$ be the column vector that is the first column of $C$. We also recall $\mathbf{I}_k$ stands for the identity matrix in dimension $k$.
Then the relation \eqref{eq:k-key} may be written in matrix terms as follows:
\begin{equation}\label{eq-Phi}\left(C \circ (\mathbf{I}_k-A) \right) W  = C_0 \circ M,
\end{equation}
using $\circ$ the Hadamard or term-by-term product (beware the two matrices we take the product of must have the same dimension). 
Using Cramer formula, one can solve \eqref{eq-Phi} and compute the expression of $W$ in terms of $A$, $M$ and the coefficients $\beta_0, \beta_1, \ldots$, taking advantage of Lemma \ref{hadamard} to simplify the resulting expression.

\begin{theorem}
It holds for every integer $i \in \K$, and complex $z$ such that $|z| < \textcolor{red}{z_c}$, denoting by $(\mathbf{I}_k-A, M)_{j}$ the matrix $\mathbf{I}_k-A$ whose $j$-th column is replaced by the column $M$:
$$w^{(j)}(z) = \frac{\prod_{q=1}^j(1-\beta_q)}{(1-\bar{\beta})^j} \cdot \frac{\det((\mathbf{I}_k-A(z),M)_{j})}{\det(\mathbf{I}_k-A(z))}\cdot $$
\end{theorem}

\begin{proof}From Equation \eqref{eq-Phi} and Cramer formula, we have 
$$w^{(j)}(z)=\frac{\det((C\circ(\mathbf{I}_k-A), C_0 \circ M)_{j})}{\det(C\circ(\mathbf{I}_k-A))},$$
provided the denominator is non null.
For the denominator, Lemma \ref{hadamard} applied to matrices $C$ and $\mathbf{I}_k-A$ gives that
$\det(C \circ (\mathbf{I}_k-A)) = \det(\mathbf{I}_k-A)$, 
so $|z|<z_c$ is indeed sufficient to guarantee the denominator is non null (it is even positive).
For the numerator, the $j$-th column  $C_j$ of the matrix $C$ satisfies:
 $$C_0=\frac{\prod_{q=1}^j(1-\beta_q)}{(1-\bar{\beta})^j} C_j \quad \text{ hence }\quad C_0 \circ M  = \frac{\prod_{q=1}^j(1-\beta_q)}{(1-\bar{\beta})^j} (C_j \circ M).$$
Therefore $\det((C\circ(\mathbf{I}_k-A), C_0 \circ M)_{j})=\frac{\prod_{q=1}^j(1-\beta_q)}{(1-\bar{\beta})^j} \det(C\circ (\mathbf{I}_k-A,M)_{j})=\frac{\prod_{q=1}^j(1-\beta_q)}{(1-\bar{\beta})^j}\det((\mathbf{I}_k-A,M)_{j})$ using Lemma \ref{hadamard} again, which concludes the proof.

\end{proof}

\subsection{Periodic selection}

We considered a model where the effect of selection is modeled through size-biasing  by the identity map at each time. In the original case $k=1$, this is indeed no loss of generality, since one can always reduce to this case. This is no more true with 2 or more environments. It is thus of interest to consider a more general case in which the strength of selection varies at each time step, as mutation probabilities and mutation fitness distribution do. Precisely, we assume that the selection strength is quantified by non-null applications $s_{0}, s_{1}, \ldots, s_{k-1}: [0,1] \to \mathbb R_+$.
We then work under the same set of assumptions than in Section \ref{periodic-sec}, namely $\beta_i \in (0,1)$, $i \in \K$ and \eqref{cond:eta}, and consider the following sequence of fitness distributions:
\begin{equation*}\label{eq:model-variable-sel} p_{n+1}(dx)  = \beta_{[n+1]} q_{[n+1]}(x) + (1-\beta_{[n+1]}) \frac{s_{[n+1]}(x) p_{n}(dx)}{\int s_{[n+1]}(x)p_n(dx)}, \qquad n \geq 0,\end{equation*}
where $[n]$ again stands for the class of $n$ in the division by $k$.
Our proof for the convergence of $(p_{kn+i}, i \in \K)$ does not apply if $s_0,\ldots, s_{k-1}$ are not all non-decreasing, since the key property of preservation of the order, Proposition \ref{lem:order} point (ii), holds no more in this case. Yet we can devise a fixed point for the map $(p_0, \ldots, p_{k-1}) \longmapsto (p_k, \ldots, p_{2k-1})$, which gives rise to the following conjecture:

\begin{conjecture}[Extension of Theorem \ref{main2} in case of a periodic mutation \emph{and} periodic selection]\label{conjecture}
   Let $s(x) =(\prod_{\ell \in \K} s_\ell(x))^{1/k}$ for short.
   Assuming the maximum of $s(x)$ on $\emph{Supp}(p_0)$ is larger than the \emph{positive} maximum of $s(x)$ on $\emph{Supp}(q_i)$ for each $i \in \K$ and reached at a unique point, say $x_0$,
   and with the matrix $A^{\texttt{s}}$
   defined\footnote{
the product $s_{[j+1]}(x) \ldots s_{[j+i]}(x)$ in the definition of $\mu^i_j(z)$ is a shorthand for $\prod_{j+1\leq \ell \leq j+i} s_{[\ell]}(x)$, and is by convention null when $i=0$; a similar remark applies for the product $s_{[i]}(x) \ldots s_{[i-j+1]}(x) = 
\prod_{0 \leq \ell \leq j-1} s_{[i-\ell]}$} as follows:
   $$
   A^{\texttt{s}}_{i,j}(z) = \mu_{j}^{[i-j]}(z) \text{ where } \mu_{j}^{i}(z) =  \int \frac{z^i s_{[j+1]}(x) \ldots s_{[j+i]}(x)}{1-z^k s(x)^k} \beta_j q_j(dx), \qquad  i,j \in \K$$
 for $z \in I_{\mu}$ defined from $\mu_{j}^{i}$ as in \eqref{eq:mu-ij}, then, for $z_c$, $\alpha$ and $U$ defined from the new matrix $A^{\texttt{s}}$ as in \eqref{U-theta}, it holds:
$$
p_{kn+i}(dx) \underset{n\to\infty}{\longrightarrow} \pi_i(dx)  \defeq \sum_{j=0}^{k-1} \frac{U_{[i-j]}}{U_{i}} \frac{z_c^j s_{[i]}(x) \ldots s_{[i-j+1]}(x) \beta_{[i-j]} q_{[i-j]}(dx)}{1- z_c^k s(x)^k}  + \alpha\frac{U_{0}}{U_{i}} \delta_{x_0},   \qquad \text{for } i \in \K,
$$
with the convergence in the same sense as in Theorem \ref{main2}.
\end{conjecture}
\begin{remark}
    In case the maximum of $\prod_{\ell \in \K} s_\ell(x)$ 
 on Supp($p_0$) is reached at multiple points, it is not clear that the mass of $p_{kn+i}(dx)$ at each of those points should converge; one possible way round this issue would then be to sum up the masses at these points.
\end{remark}

\appendix

\section{Appendix}\label{sec-appendix}

In this Appendix, and apart from the very last proof, the symbols $A, B, \ldots$ will denote generic matrices.
For the ease of reference, we start the Appendix by recalling the Perron Theorem together with a set of properties that we shall need in this article. The reader interested in the rich history of the theorem may indulge himself in the reading of \cite{hawkins2008continued} (or of the shorter piece \cite{maccluer2000many}).

\begin{theorem}[Perron, \cite{Perron1907}]\label{thm:Perron}
    Let $A=(A_{ij})_{0\leq i,j < k}\in \mathcal{M}_k(\mathbb{R}_+^{\star})$ be a matrix with positive coefficients. 
    \begin{enumerate}
        \item There is a positive real number $r$, called the Perron eigenvalue, such that $r$ is an eigenvalue of $A$ and any other eigenvalue $\lambda \in \mathbb C$ of $A$ is such that $|\lambda|< r$.
        \item The Perron eigenvalue is a simple root of the characteristic polynomial of $A$, and the associated eigenspace is one-dimensional.
        \item There is an eigenvector $V$ of $A$ associated with the eigenvalues $r$ whose coordinates are all positive.
        \item There are no other eigenvectors (except multiples of $V$) of $A$  
        whose coordinates are all non-negative.
    \end{enumerate}
\end{theorem}

\begin{proposition}\label{prop-Perron}
    \begin{description}
     
        \item[$(i)$] (Collatz-Wielandt) 
        The Perron eigenvalue $r$ of $A=(A_{ij})_{0\leq i,j < k}\in\mathcal{M}_k(\mathbb{R}_+^{\star})$ satisfies
        $$r= \max_{X} \min_{i: X_i \neq 0} \frac{(AX)_i}{X_i}$$
        where the max is taken over those vectors $X$ whose coefficients are all non-negative and not all null. 
        \item[$(ii)$] (Row sums bounds)        
        It also satisfies $\min_{i} \sum_{j} A_{i,j}  \leq r \leq \max_{i} \sum_{j} A_{i,j}$.
        \item[$(iii)$] (Continuity) The map $A \mapsto r(A)$ is continuous from $\mathcal{M}_k(\mathbb{R}_+^{\star})$ to $\mathbb{R}_+^{\star}$.
        \item[$(iv)$] If $A,B\in \mathcal{M}_k(\mathbb{R}_+^{\star})$ are such that $0<A_{i,j}< B_{i,j}$ holds for each $i,j\in\K$ then their respective Perron eigenvalues $r(A)$ and $r(B)$ satisfy $$r(A)<r(B).$$
    \end{description}
\end{proposition}

\begin{proof}
The proof of points $(i)$ to $(iii)$ of Proposition \ref{prop-Perron} can be found in \cite{Meyer2023matrix}, respectively at p. 666, p. 668 (8.2.7) and p. 497 (7.1.3). We therefore only detail here the proof of point $(iv)$ of Proposition \ref{prop-Perron}. Let $Y$ be the eigenvector of $B$ associated to the Perron eigenvalue $\rho(B)$. Then $Y$ has positive entries and for each $i \in \K$, $$\rho(B)Y_i=\sum_jB_{i,j}Y_j<\sum_{j}A_{i,j}Y_j=(AY)_i.$$ 
Collatz-Wielandt formula from point $(i)$ then gives
$$\rho(A)= \max_{X} \min_{i: X_i \neq 0} \frac{(AX)_i}{X_i} \geq \min_{i: Y_i \neq 0} \frac{(AY)_i}{Y_i} > \rho(B).$$
\end{proof}

\begin{proposition}[Adjugate matrix]
\label{adjugate}
Let $B \in \mathcal{M}_k(\mathbb{R})$. The adjugate matrix $\emph{adj}(B)$ of $B$, defined by $\emph{adj}(B) = (-1)^{i+j} M_{ji},$ where the $(i,j)$-minor $M_{ij}$ is the determinant of the matrix obtained from $B$ by deleting its $i$-th line and $j$-th column, satisfies:
\begin{equation}
\label{adj-inv}
B \;\emph{adj}(B) = \emph{adj}(B) \; B = \emph{det}(B) \mathbf{I}_k.
\end{equation}
As a consequence, one has the following trichotomy:
 \begin{description}
\item In case $\emph{rank}(B)=k$, $\emph{rank}(\emph{adj}(B))=k$: $B$ is invertible. 
\item In case $\emph{rank}(B)=k-1$, $\emph{rank}(\emph{adj}(B))=1$.
\item In case $\emph{rank}(B)\leq k-2$, $\emph{rank}(\emph{adj}(B))=0$ : $B$ is the null matrix.
\end{description}
\end{proposition}

\begin{proof}
The relation \eqref{adj-inv} is classical (by e.g. 6.2.7 in \cite{Meyer2023matrix}). We give a short justification for the only non-direct point in the ensuing trichotomy for the ease of reference: if rank$(B)=k-1$, from \eqref{adj-inv}, the image of adj$(B)$ is included in the one-dimensional kernel of $A$, hence rank$(B) \leq 1$; but one of the $(i,j)$-minors $M_{ij}$ of $B$ has to be non null by the rank assumption, hence adj$(B)$ is not the null matrix and rank$(B)$ is $1$. 

\end{proof}

\begin{lemma}
\label{hadamard}
Let $A = (A_{i,j})_{0 \leq i,j <k} \in\mathcal{M}_k(\mathbb{R})$ and $\beta_0,...,\beta_{k-2}$ be a collection of non-null real numbers. Define another matrix $B =  (B_{i,j})_{0 \leq i,j <k} \in\mathcal{M}_k(\mathbb{R})$ as follows: $B_{ii}=1$ for $0 \leq i <k$,  $B_{ij}=\prod_{k=j}^{i-1}\beta_k$ for $0 \leq j<i<k$ and $B_{ij}=1/B_{ji}$  for $0 \leq i <j <k$. The determinant of the matrix $B \circ A$, defined as the Hadamard (or term-by-term) product between $A$ and $B$, does not depend on the quantities $\beta_i$, i.e. 
$$\det(B \circ A)=\det(A).$$
\end{lemma}

\begin{proof}
This lemma can be proved by induction on the dimension $k$ of the matrix. First the result is true for $k=2$ since $A_{00}A_{11}-\beta_0 A_{10} \times (1/\beta_0)  A_{01}=A_{00}A_{11}-A_{10}A_{01}$.
Now we detail the induction.
For $M \in \mathcal{M}_k(\mathbb{R})$ a matrix, let $M^{(ij)} \in \mathcal{M}_{k-1}(\mathbb{R})$  be the matrix obtained from $M$ by erasing its $i$-th line and $j$ column.
Expanding the determinant of $C:=B \circ A$ according to its first column and using the definition of $C$ one has:
\begin{align*}
\det(C)& =\sum_{i=0}^{k-1} (-1)^{i}A_{i0}\left(\prod_{k=0}^{i-1}\beta_k\right) \; \det(C^{(i0)}) \\
& =\sum_{i=0}^{k-1} (-1)^{i} A_{i0} \; \det(A^{(i0)})\\
& =\det(A).
\end{align*}
To go from the first to the second line, observe that the matrix $C^{(i0)}$ is such that multiplying the $i-1$ first lines respectively by $\beta_0$, $\beta_1$, ..., $\beta_{i-2}$  gives a matrix of size $k-1$ with the same form as $C$, but with a reduced size, hence the second line by the induction assumption.
\end{proof}

\begin{lemma}
\label{invariance-rotation}
Let $A = (A_{i,j})_{0 \leq i,j <k} \in\mathcal{M}_k(\mathbb{R})$, and define $s(A)= (A_{[i+1],[j+1]})_{0 \leq i,j<k }$  where $[i]$ stands for $i$ modulo $k$. Then 
$$\det(s(A))=\det(A).$$
\end{lemma}

\begin{proof}
Observe that $s(A)$ may be built from $A$ by pushing the first column of $A$ after the last column of $A$, and then, pushing the first line of the newly formed matrix after its last line ; in terms of determinants, each of these push is equivalent to exchanging $k-1$ columns/lines, hence modifies the value of the determinant by $(-1)^{k-1}$; the two pushes hence leave the determinant invariant.
\end{proof}

Note that if $A(z)$ is the matrix associated with the mutation environment $(\beta_i q_{i})_{i \in \K}$, $s(A(z))$ is the matrix associated with the mutation environment $(\beta_{[i+1]} q_{[i+1]})_{i \in \K}$. This entails that $\text{det}(\mathbf{I}_k-A(z))$ and $\rho(A(z))$ are invariant under cyclic permutations of the environments as discussed in Section \ref{periodic-sec}.

\begin{lemma}
\label{1-not-eigenvalue}
    Let $B\in\mathcal{M}_k(\mathbb{R_+})$ such that for any $i\in \K $, 
    \begin{equation}
    \label{generic-ass}
        B_{i,i} > B_{i, [i+1]} > \ldots > B_{i,[i+k-1]}  > \max\{B_{i,i}-1,0\}
    \end{equation} 
    Then there is no vector $U \in \mathbb R^k$ admitting both positive and negative coordinates for which
    $$BU=\lambda U \text {  for a real number } \lambda \geq 1.$$ 
\end{lemma}
 
\begin{proof} Wlog, one can assume that $\sum_{j\in\K}U_j\geq 0$. Wlog again, up to rotation of the matrix at the first index where the coordinate of $U$ is negative, one can assume that $U_0<0$ \footnote{formally, we replace $B$ by $C= (B_{[r_0+i],[r_0+j]})_{i,j\in\K}$, where $r_0$ is the first index such that $U_{r_0}>0$}. For such a vector $U$, the minimum of the sequence of partial sums $(U_0, U_{0}+U_{1},...,\sum_{j\in\K} U_j)$ is negative. Then, up to rotation of the matrix at the first index such that the partial sum is miminal (as in the so-called cyclic lemma), one can assume that $\sum_{j=1}^\ell U_j \geq 0$ for each $\ell\in\K$, and $U_{k-1}< 0$ \footnote{this time, we replace $C$ by $D=(C_{[r_1+i],[r_1+j]})_{i,j\in\K}$, where $r_1$ is the first index such that  $\sum_{0 \leq i \leq j} U_i$ is minimal}, which also implies $\sum_{i=0}^{k-2} U_i > 0$. 

We now have that $\lambda U_{k-1}=\sum_{j=0}^{k-1}B_{k-1,j} U_{j}$ for a new matrix (abusively) denoted $B$ again, that still satisfies \eqref{generic-ass} hence:
\begin{align*}(\lambda-B_{k-1,k-1})U_{k-1}&=\sum_{j=0}^{k-2} B_{k-1,j} U_j\\
&=B_{k-1,k-2}\sum_{j=0}^{k-2}U_j+(B_{k-1,k-3}-B_{k-1, k-2})\sum_{j=0}^{k-3} U_j+...+(B_{k-1,0}-B_{k-1,1}) U_0 \\ 
&\geq B_{k-1,k-2} \sum_{j=0}^{k-2}U_j\end{align*}
using for the inequality our assumption that $B_{k-1,0}>B_{k-1,1}> \ldots > B_{k-1,k-2}$ together with the fact that $\sum_{j=0}^{\ell} U_j \geq 0$ for each $\ell$.
Dividing by the positive quantity $B_{k-1,k-2}$, and using the assumption that $\lambda\geq 1$ together with $B_{k-1,k-2} >B_{k-1,k-1}-1$, we finally get:
$$\sum_{i=0}^{k-2}U_i \leq \frac{B_{k-1,k-1}-\lambda}{B_{k-1,k-2}}(-U_{k-1}) \leq
\frac{B_{k-1,k-1}-1}{B_{k-1,k-2}}(-U_{k-1}) < -U_{k-1},$$
that is, $\sum_{i=0}^{k-1}U_i<0$, which is in contradiction with our first assumption about that sum.
\end{proof}

We conclude with the proof of Proposition \ref{at-most-one}. Recall the aim is to prove that, for $0 \leq z \leq 1/\eta_{q}$, $A(z)$ has at most one real eigenvalue that is $\geq 1$.

\begin{proof}[Proof of Proposition \ref{at-most-one}]
Let $z>0$. From the definition
$A_{i,j}(z) = \int \frac{(zx)^{[i-j]}}{1-(zx)^k}\beta_{j} q_{j}(dx)$
and the fact that $\beta_j q_j(dx)$ is a sub-probability measure, we get, for $\ell \in \{0, \ldots, k-2\}$, using that $0 \leq zx \leq 1$,
\begin{align*} 
A_{[i+\ell],i}(z) - A_{[i+\ell+1],i}(z) & = \int \frac{(zx)^{\ell} (1-zx)}{1-(zx)^k}\beta_{i} q_{i}(dx) >0 \\
A_{[i+k-1],i}(z) &= \int \frac{(zx)^{k-1}}{1-(zx)^k}\beta_{i} q_{i}(dx) >0 \\
A_{i,i}(z)- A_{[i+k-1],i}(z) & =\int \frac{1-(zx)^{k-1}}{1-(zx)^k}\beta_{i} q_{i}(dx) <1
 \end{align*}
 This shows the transpose $B(z)$ of the $A(z)$ matrix satifies the assumption \eqref{generic-ass} of  Lemma \ref{1-not-eigenvalue}. As a consequence, there is no real-valued eigenvector of $B$ with both positive and negative coordinates 
 associated with a real eigenvalue $\lambda \geq 1$. But: i) any eigenvector of a real-valued matrix associated with a real eigenvalue has real coordinates ii) any eigenvector distinct from the Perron eigenvector admits both positive and negative coordinates (this is the contraposition of Theorem \ref{prop-Perron} point 4). Hence, the matrix $B(z)$ satisfies the statement of Proposition \ref{at-most-one}, and the same is then true of its transpose matrix $A(z)$.
\end{proof}

\bibliographystyle{abbrv}
\bibliography{Biblio}

\end{document}